\documentclass[12pt]{article}
\usepackage[letterpaper, left=1.5in, right=1in,top=1in,bottom=1in]{geometry}
\usepackage{authblk}
\usepackage{tikz}
\usepackage{chngcntr}
\counterwithin{figure}{section}
\usetikzlibrary{arrows}
\usepackage{placeins}	
\usepackage{amsfonts}
\usepackage{graphicx}
\usepackage{standalone}
\usepackage{amsmath, amsthm, amssymb}
\usetikzlibrary{decorations.markings}
\tikzset{->-/.style={decoration={
markings,
mark=at position .55 with {\arrow{>}}},postaction={decorate}}}
\tikzset{-<-/.style={decoration={
markings,
mark=at position .55 with {\arrow{<}}},postaction={decorate}}}
\usepackage[capitalize, noabbrev]{cleveref}
\newtheorem{thrm}{Theorem}[section]

\newtheorem{lem}[thrm]{Lemma}
\newtheorem{cor}[thrm]{Corollary}
\newtheorem{defn}[thrm]{Definition}

\newenvironment{pf}           {\noindent{\bf Proof:} }%
                                {\null\hfill$\Box$\par\medskip}
                           
\usepackage[nottoc]{tocbibind}

\newcommand{\rell}[2]{\mbox{SCNRel$(#1,#2)$}}
\newcommand{\bfam}[2]{\overleftrightarrow{\mathcal{D}}_{#1,#2}}
\newcommand{\nbfam}[2]{\overrightarrow{\mathcal{D}}_{#1,#2}}

\begin{document}

\title{Existence of Optimally-Greatest Digraphs for Strongly Connected Node Reliability}
\author[1]{Danielle Cox}
\author[2]{Kyle MacKeigan}
\author[3]{Emily Wright}
\affil[1]{Department of Mathematics, Mount Saint Vincent University}
\affil[ ]{Email: \textit {Danielle.cox@msvu.ca}}
\affil[2]{Department of Mathematics, Mount Saint Vincent University}
\affil[ ]{Email: \textit{Kyle.m.mackeigan@gmail.com}}
\affil[3]{Department of Mathematics, Concordia University}
\affil[ ]{Email: \textit{emily.wright@concordia.ca}}
\maketitle

\begin{abstract}

In this paper, we introduce a new model to study network reliability with node failures. This model, strongly connected node reliability, is the directed variant of node reliability and measures the probability that the operational vertices induce a subdigraph that is strongly connected. If we are restricted to directed graphs with $n$ vertices and $n+1\leq m\leq 2n-3$ or $m=2n$ arcs, an optimally-greatest digraph does not exist. Furthermore, we study optimally-greatest directed circulant graphs when the vertices operate with probability $p$ near zero and near one. 

In particular, we show that the graph $\Gamma\left(\mathbb{Z}_n,\{1,-1\}\right)$ is optimally-greatest for values of $p$ near zero. Then, we determine that the graph $\Gamma\left(\mathbb{Z}_{n},\{1,\frac{n+2}{2}\}\right)$ is optimally-greatest for values of $p$ near one when $n$ is even. Next, we show that the graph $\Gamma\left(\mathbb{Z}_{n},\{1,2(3^{-1})\}\right)$ is optimally-greatest for values of $p$ near one when $n$ is odd and not divisible by three and that $\Gamma\left(\mathbb{Z}_{n},\{1,3(2^{-1})\}\right)$ is optimally-greatest for values of $p$ near one when $n$ is odd and divisible by three. We conclude with a discussion of open problems.
\end{abstract}

\section{Introduction}

Network reliability is a well studied area of graph theory and there are a variety of models used to study this problem. The most general model is as follows. Let $G$ be a graph and let $A$ be a set of elements that are independently operational with probability $p\in [0,1]$. Then, for some graph property $\mathcal{P}$, let $Rel_{\mathcal{P}}(G,p)$ denote the probability that the operational elements of $A$ induce a subgraph of $G$ that has property $\mathcal{P}$. This probability can be explicitly computed as

\[ Rel_{\mathcal{P}}(G,p)=\sum_{i=0}^{|A|}F_i(G)p^{|A|-i}(1-p)^i\]

where the coefficient $F_i(G)$ counts the number of ways that $i$ elements of $A$ can be non-operational and the resulting induced subgraphs have property $\mathcal{P}$. This is the $F$-form of the reliability polynomial. Note we write $F_i$ when $G$ is clear from context. Similarly, we have the $N$-form of the reliability polynomial, 
\[ Rel_{\mathcal{P}}(G,p)=\sum_{i=0}^{|A|}N_i(G)p^{i}(1-p)^{|A|-i}\]

where $N_i(G)$ represents the number of ways to have $i$ elements of $A$ operational such that an induced graph with property $\mathcal{P}$ results. We write $N_i$ when $G$ is clear from context. It is easy to see that for $0\leq i\leq |A|$, $N_i=F_{|A|-i}$.

Since this summation is a polynomial in $p$, it is called a \emph{reliability polynomial}. For a recent survey on reliability polynomials see \cite{browncoxhighway}, and for a background on network reliability see \cite{colbourn}. We now provide a brief summary on several well studied models of reliability.

The \emph{all-terminal reliability} of a graph $G$ is the model where $A=E(G)$, the edge set, and $\mathcal{P}$ is the property that at least a spanning tree is operational (see \cite{colbourn}, for example). This model ensures that all the vertices of $G$ can communicate with one another via the operational edges. The \emph{node reliability} of a  graph $G$ is the model where $A=V(G)$, the vertex set, and $\mathcal{P}$ is the property that a connected subgraph is operational (see \cite{colbourn,stivaros,sutner1991complexity} for example).  The directed version of all-terminal reliability is called \emph{strongly connected reliability} and was introduced in \cite{brown2005strongly} which studied the existence of optimally-greatest digraphs for this model. Other areas of research for this model were regarding the analytic properties \cite{browncoxroots,browndilcher,coxphdthesis}.

In this paper we introduce the directed version of node reliability, the \emph{strongly connected node reliability} of a digraph, $G$, denoted $\rell{G}{p}$.  We will provide preliminary results regarding the polynomial, focus on the optimality of directed circulants and highlight differences between the undirected and directed versions of node reliability.

\subsubsection{Strongly Connected Node Reliability}
The network reliability model \emph{node reliability} was introduced by Sutner et al. in 1991 \cite{sutner1991complexity}, and in recent years has applications to the study of to social media networks. If we let $G$ be a graph on $n$ vertices and $m$ edges and assume that the vertices of $G$ operate independently with probability $p\in[0,1]$, the node reliability of $G$, NRel$(G,p)$ is the probability that the subgraph induced by the operational vertices is connected. 

Many real-world networks can be modelled as a digraph where vertices represent hubs where information is processed, then transmitted throughout a network. For example, with respect to social media, the operational vertices model users being logged onto the social media platform, thus able to direct message each other and get immediate responses. Relationships between users on platforms like LinkedIn and Facebook are two way, two users must agree to follow one another. For platforms like Twitter, TikTok and Instagram, the relationship between users may only be one way - person A can follow person B, but person B may choose to not follow person A. We begin with the following formal definition of strongly connected node reliability.

\begin{defn}
Let $G$ be a strongly connected digraph on $n$ vertices and $m$ arcs. Assume that vertices operate independently with probability $p\in[0,1].$ The \textbf{strongly connected node reliability} of a digraph $G$, denoted $\rell{G}{p}$, is the probability that the subdigraph induced by the operational vertices is strongly connected. 
\end{defn}

For example, for any directed cycle $C^{\rightarrow}_n$, $n\geq 2$ we require either exactly one node to be operational or all nodes. Thus, $\rell{C^{\rightarrow}_n}{p}=np(1-p)^{n-1}+p^n$. In general, considering the $F$-form of the strongly connected node reliability polynomial $F_0=1$, $F_1=n$ and $F_2=b$, where $b$ is the number of \emph{bundles}, which are pairs of antiparallel arcs that exist between two distinct vertices. In the next section we study the existence of optimally-greatest digraphs, which will require a study of the coefficients of the $F$-form of the strongly connected node reliability polynomial.

\section{Optimally-Greatest Digraphs, $m\leq 2n-1$ Arcs}

In this section we determine the existence of digraphs on $n$ vertices and $m$ edges whose reliability equals or exceeds that of any other graph of the same order and size  for all $p\in [0,1]$. 

\begin{defn}
Let $\mathcal{F}$ be a family of graphs. We call $G\in\mathcal{F}$ an $\mathcal{F}$-graph, and say that $G$ is an optimally-greatest $\mathcal{F}$-graph if for any other $\mathcal{F}$-graph $H$ we have that $\rell{G}{p}\geq \rell{H}{p}$ for all $p\in [0,1]$. Optimally-least is similarly defined.
\end{defn}

This is a global notion of optimality. If we know that the nodes are highly reliable, or highly unreliable, we may be interested in a more local notion of optimality. For $p_{0} = 0$ or $1$ we will also talk about $G$ being optimally-greatest $\mathcal{F}$-graph sufficiently close to $p_{0}$ when there is an $\varepsilon > 0 $ such that $\rell{G}{p}\geq \rell{H}{p}$ for all $p\in [0,1] \cap (p_{0}-\varepsilon,p_{0}+\varepsilon)$ and for all other $\mathcal{F}$-graphs $H$.

In this work we will only be considering digraphs that are strongly connected, since otherwise, if all the vertices are operational, the strongly connected node reliability would be identically 0. The following observations from \cite{browncoxoptimal} will be useful.

\begin{lem}\cite{browncoxoptimal} 
\label{Lemma: N compare}
Let $G,H\in \mathcal{F}_{n,m}$. Consider the $N$-form of the reliability polynomial: 
\[ \rell{G}{p}=\sum_{i=1}^{n}N_i(G)p^{i}(1-p)^{n-i}  \hspace{30pt}\rell{H}{p}=\sum_{i=1}^{n}N_i(H)p^{i}(1-p)^{n-i} .\]
If there is a $k$ such that for all $i<k$, $N_i(G)=N_i(H)$, and for $i=k$, $N_k(G)>N_k(H)$, then $\rell{G}{p}>\rell{H}{p}$ for values of $p$ sufficiently close to 0.
\end{lem}

\begin{lem}\cite{browncoxoptimal} 
\label{Lemma: F compare}
Let $G,H\in \mathcal{F}_{n,m}$. Consider the $F$-form of the reliability polynomial: 
\[ \rell{G}{p}=\sum_{i=0}^{n-1}F_i(G)p^{n-i}(1-p)^i   \hspace{30pt}\rell{H}{p}=\sum_{i=0}^{n-1}F_i(H)p^{n-i}(1-p)^i .\]
If there is a $k$ such that for all $i<k$, $F_i(G)=F_i(H)$, and for $i=k$, $F_k(G)>F_k(H)$, then $\rell{G}{p}>\rell{H}{p}$ for values of $p$ sufficiently close to 1.
\end{lem}

Let $\bfam{n}{m}$ represent the family of all strongly connected digraphs on $n$ vertices and $m$ arcs, allowing bundles, and $\nbfam{n}{m}$ the family of digraphs of size $m$ and order $n$ that do not have bundles. For these two families, given a fixed $n\geq 2$, we want to know for what values of $m$ does there exist a digraph $D$ that is optimally-greatest. 

Clearly, if $m\leq n-1$, then a strongly connected digraph does not exist. If $n=m$, then the only possibly digraph is the directed cycle, and hence is optimally-greatest for both $\bfam{n}{n}$ and $\nbfam{n}{n}$

\begin{lem}
\label{sparsept1}
Let $n \geq 5$ and $m=2n-k$, where $3 \leq k \leq n-1$ (therefore $n+1 \leq m \leq 2n-3$). If  an optimally-greatest digraph, $G$ exists for $\bfam{n}{m}$ then $U(G)$ is unicyclic and $G$  has $n-k$ bundles. If $m=2n-2$, then $G$ has $n-1$ bundles.
\end{lem}

\begin{pf}
Let $G$ be an optimally-greatest digraph $\bfam{n}{m}$. Thus it is optimally-greatest for values of $p$ near zero, therefore by Lemma~\ref{Lemma: N compare}, it must have the maximum number of bundles possible (maximize $N_2$). 
Suppose that $G$ has $n-k+\ell$, $\ell \geq 0$ bundles and $k-2\ell$ non-bundle arcs. It follows that $U(G)$, the underlying undirected graph of $G$ has $n$ vertices and $n-k+\ell+k-2\ell=n-\ell$ edges.

If $\ell \geq 2$ then $U(G)$ is not connected, hence a contradiction as $G$ is strongly connected.

If $\ell =1$ then $U(G)$ is a tree, hence all edges correspond to a bundle in $G$ and thus $G$ has $n-1$ bundles and a total of $2n-2$ arcs. If $G$ has $2n-2$ arcs, it can have at most $n-1$ bundles, hence $U(G)$ must be a tree if $G$ is optimally-greatest.

If $\ell=0$ then $G$ has $n-k$ bundles and $k$ non-bundle arcs. This means that $U(G)$ has $n$ edges, and hence is unicyclic. \end{pf}

We will now use this result to prove that for $\bfam{n}{m}$ no optimally-greatest digraph exists when $n\geq 5$ and $n+1 \leq m \leq 2n-3$.

\begin{thrm}
\label{sparsept2}
If $n \geq 5$ and $m=2n-k$, where $3 \leq k \leq n-1$, then an optimally-greatest digraph does not exist.
\end{thrm}

\begin{pf}
We will proceed by showing that if an optimally-greatest digraph exists for $n\geq 5$ and $n+1\leq m \leq 2n-3$, then it is unique by considering optimality for values of $p$ near 0. We will then find another digraph of order $n$ and $m$ that is more optimal for values of $p$ near 1.

Supposed $G$ is optimally-greatest. From Lemma~\ref{sparsept1}, we know $U(G)$ is unicyclic and $G$ has $n-k$ bundles.

All potential optimally-greatest digraphs have $N_2=n-k$, so since $G$ is optimally-greatest then it has the most number of strongly connected subdigraphs of order 3, that is has the largest possible value for $N_3$ over all digraphs in $\bfam{n}{m}$. These subdigraphs are such that three vertices, $x,y,z$, induce a path with bundles between $x$ and $y$, and $y$ and $z$, or they lie on a cycle of length 3 (which may or may not contain bundles between pairs of $x$,$y$,$z$). 
 
Let $G_k$ be the graph that is a directed cycle of length $k$, with each of the remaining $n-k$ vertices (call this set of vertices $V'$) all incident to the same vertex, $v_0$ in the cycle via a bundle. We will use induction to show that this graph has the largest value of $N_3$ over all digraphs in $\bfam{n}{2n-k}$ that have $n-k$ bundles. It is the case that $G_k$ has the largest $N_3$ for $n=5$. For illustration, the digraph $G_8$ when $n=10$ is given in Figure \ref{Figure: Gk}

\definecolor{sqsqsq}{rgb}{0.12549019607843137,0.12549019607843137,0.12549019607843137}
\begin{figure}[h!]
\centering
\begin{tikzpicture}[line cap=round,line join=round,>=triangle 45,x=1.0cm,y=1.0cm]
\draw [shift={(4.516,4.)},->-]  plot[domain=2.5584769995314325:3.7247083076481537,variable=\t]({1.*1.8161101288192858*cos(\t r)+0.*1.8161101288192858*sin(\t r)},{0.*1.8161101288192858*cos(\t r)+1.*1.8161101288192858*sin(\t r)});
\draw [shift={(1.5981249999999996,4.)},->-]  plot[domain=-0.6196166013789108:0.6196166013789107,variable=\t]({1.*1.72199114853271*cos(\t r)+0.*1.72199114853271*sin(\t r)},{0.*1.72199114853271*cos(\t r)+1.*1.72199114853271*sin(\t r)});
\draw [shift={(4.,1.0375)},->-]  plot[domain=1.0995346517388256:2.0420580018509673,variable=\t]({1.*2.2025908040305597*cos(\t r)+0.*2.2025908040305597*sin(\t r)},{0.*2.2025908040305597*cos(\t r)+1.*2.2025908040305597*sin(\t r)});
\draw [shift={(4.,4.125789473684211)},->-]  plot[domain=3.9860949604923372:5.438683000277042,variable=\t]({1.*1.505789473684211*cos(\t r)+0.*1.505789473684211*sin(\t r)},{0.*1.505789473684211*cos(\t r)+1.*1.505789473684211*sin(\t r)});
\draw [->-] (0.,0.)-- (0.,1.5);
\draw [->-] (0.,1.5)-- (0.,3.);
\draw [->-] (0.,3.)-- (1.5,3.);
\draw [->-] (1.5,3.)-- (3.,3.);
\draw [->-] (3.,3.)-- (3,1.5);
\draw [->-] (3,1.5)-- (3.,0.);
\draw [->-] (3.,0.)-- (1.5,0.);
\draw [->-] (1.5,0.)-- (0.,0.);
\begin{scriptsize}
\draw [fill=sqsqsq] (0.,0.) circle (2.5pt);
\draw [fill=sqsqsq] (0.,3.) circle (2.5pt);
\draw [fill=sqsqsq] (3.,3.) circle (2.5pt);
\draw [fill=sqsqsq] (3.,0.) circle (2.5pt);
\draw [fill=sqsqsq] (3.,5.) circle (2.5pt);
\draw [fill=sqsqsq] (5.,3.) circle (2.5pt);
\draw [fill=sqsqsq] (0.,1.5) circle (2.5pt);
\draw [fill=sqsqsq] (1.5,3.) circle (2.5pt);
\draw [fill=sqsqsq] (3.0,1.5) circle (2.5pt);
\draw [fill=sqsqsq] (1.5,0.) circle (2.5pt);
\end{scriptsize}
\end{tikzpicture}
\caption{Digraph $G_8$}
\label{Figure: Gk}
\end{figure}
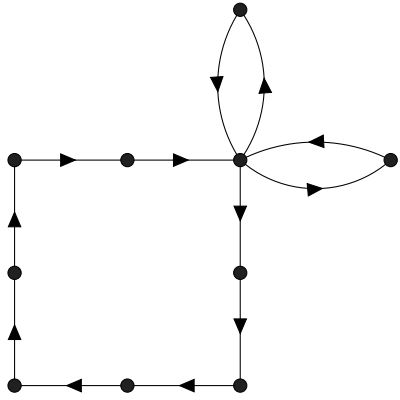

Suppose that $n\geq 6$. In $G_k$ let $x\in V'$ be a vertex of $G_k$. There are $n-(k+1)$ strongly connected subdigraphs of order 3 that contain $x$, namely, the digraph induced by the distinct vertices $x,v_0,z$ for $z\in V'$ and ${n-(k+1)\choose 2}$ without $x$, which consist of $v_0$ and two vertices from $V'\setminus {x}$. It should be noted if $k=3$ there is an additional subdigraph of order 3.
 
Let $H$ be another digraph with $n-k$ bundles in $\bfam{n}{2n-k}$, for $3\leq k\leq n-1$. First assume $U(H)$ is such that there exists some vertex $y$ of degree one, hence in $H$ adjacent to one other vertex via a bundle. The number of strongly connected subdigraphs of order 3 of $H$ that contain $y$ is at most $n-(k+1)$ as the sole neighbour of $y$ must be operational and there are $n-k-1$ bundles in the digraph not including those incident to $y$. Note, equality is achieved only when $H=G_k$. Consider $H-y$, this digraph has $2n-k-2=2(n-1)-k$ arcs, $3\leq k\leq n-1$ and $n-1$ vertices. If $3\leq k \leq (n-1)-1=n-2$ by the induction hypothesis $H-y$ has strictly less strongly connected subdigraphs of order 3 than $G_{k-1}$. Suppose that $k=n-1$. The induction hypothesis does not apply in this case, but $H-y$ has $n-1$ vertices and $m=2n-k-2=2n-(n-1)-2=n-1$ arcs, and must be a directed cycle else $H$ is not strongly connected, hence $H=G_k$ is the unique optimally-greatest digraph in this case. Thus $H$ has at most $n-k-1+N_3(G_{k-1})\leq N_3(G_k)$, with equality when $H$ is isomorphic to $G_k$.

Now suppose $U(H)$ is has no vertices of degree 1. If $H$ is also an optimally-greatest digraph by Lemma~\ref{sparsept1} we know $H$is unicyclic and has $n-k$ bundles thus $U(H)$ must have $2n-k-(n-k)=n$ edges, so must be a cycle of order $n\geq 6$. This means that for $H$ all the strongly connected subdigraphs of order 3 are paths of length 3 with bundles between adjacent vertices, thus $N_3(H)\leq n-k-2<N_3(G_k)=n-k-1+{n-(k+1) \choose 2}$. Hence $G_k$ is the unique digraph with the largest value of $N_3$.
 
Let $H_k$ be a family of digraphs also on $n$ vertices and $m=2n-k$ arcs, $3 \leq k \leq n-1$, where $U(H_k)$ is a theta graph.
Let $H_k$ have two vertices $x$ and $y$, and let there be $n-k+1$ paths of length two from $x$ to $y$. Therefore, we have $n-(n-k+3)=k-3$ vertices and $2n-k-(2n-2k+2)=k-2$ arcs remaining. 
To ensure that $H_k$ is strongly connected the remaining $k-2$ arcs  are used so that the remaining $k-3$ vertices lay in a path directed from $y$ to $x$. This means $H_k$ has $n-k+1$ ways to remove one vertex and still induce a strongly connected subdigraph. For illustration, the digraph $H_8$ when $n=10$ is given in Figure \ref{Figure: Hk}.

\begin{figure}[h!]
\centering
\begin{tikzpicture}[line cap=round,line join=round,>=triangle 45,x=1.0cm,y=1.0cm]
\draw [->-] (6.,0.)-- (5.,0.);
\draw [->-] (5.,0.)-- (4.,0.);
\draw [->-] (4.,0.)-- (3.,0.);
\draw [->-] (3.,0.)-- (2.,0.);
\draw [->-] (2.,0.)-- (1.,0.);
\draw [->-] (1.,0.)-- (0.,0.);
\draw [->-] (0.,0.)-- (3.,1.);
\draw [->-] (3.,1.)-- (6.,0.);
\draw [->-] (0.,0.)-- (3.,2.);
\draw [->-] (3.,2.)-- (6.,0.);
\draw [->-] (0.,0.)-- (3.,3.);
\draw [->-] (3.,3.)-- (6.,0.);
\begin{scriptsize}
\draw [fill=sqsqsq] (0.,0.) circle (2.5pt);
\draw [fill=sqsqsq] (6.,0.) circle (2.5pt);
\draw [fill=sqsqsq] (3.,1.) circle (2.5pt);
\draw [fill=sqsqsq] (3.,2.) circle (2.5pt);
\draw [fill=sqsqsq] (3.,3.) circle (2.5pt);
\draw [fill=sqsqsq] (1.,0.) circle (2.5pt);
\draw [fill=sqsqsq] (2.,0.) circle (2.5pt);
\draw [fill=sqsqsq] (3.,0.) circle (2.5pt);
\draw [fill=sqsqsq] (4.,0.) circle (2.5pt);
\draw [fill=sqsqsq] (5.,0.) circle (2.5pt);
\end{scriptsize}
\end{tikzpicture}
\caption{Digraph $H_8$}
\label{Figure: Hk}
\end{figure}
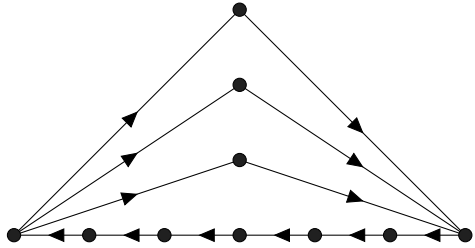

Since $G_k$ only has $n-k$ ways to remove one vertex and still induce a strongly connected subdigraph, thus by Lemma~\ref{Lemma: F compare} $H_k$ is optimally greater for values of $p$ near 1, thus, no optimally-greatest digraph exists.
\end{pf}

Let us now consider digraphs with $2n-2$ arcs, and show to contrast Theorem~\ref{sparsept2}, that an optimally-greatest digraph does exist. Let $S_n$ be the digraph whose underlying graph is a star on $n$ vertices and each edge is replaced with a bundle.

\begin{thrm}\label{thrm:2n-2}
Let $G$ be a strongly connected digraph of order $n\geq 3$ and $m=2n-2$. Then $\rell{S_n}{p}\geq \rell{G}{p}$ for all $p\in[0,1]$.
\end{thrm}

\begin{pf}
Let $G$ be a strongly connected digraph of order $n\geq 3$ and $m=2n-2$. An optimally-greatest graph on $n$ vertices and $2n-2$ arcs must have the largest number of bundles by Lemma~\ref{Lemma: F compare}, thus must be a digraph whose underlying graph is a tree and each edge is replaced by a bundle.

For $S_n$ and $G\in \bfam{n}{2n-2}$, consider $\rell{S_n}{p}-\rell{G}{p}$. We will show that this different is positive for $G$ not isomorphic to $S_n$. We will proceed by induction on $n$. It is straightforward to check that for $n=3$, $S_3$ is optimally-greatest. Suppose $n\geq 4$.

Let $v$ be a leaf of $U(S_n)$. Then \[ \rell{S_n}{p}=(1-p)\rell{S_{n-1}}{p}+p((1-p)^{n-1}+p)\] as either $v$ operational or not. If $v$ is operational then it is either the only vertex operational, or the central vertex of $U(S_n)$ must be operational. Otherwise $v$ is non-operational and we need $S_{n-1}$ to have a strongly connected subdigraph operational.

We first show that $G$ must have a cut vertex. If not, then every vertex must have in-degree (notated by deg$_+$(v)) (and out-degree) at least 2, which means \[ m=2n-2=\sum_{v\in V(D)}\mbox{deg}_+(v) \geq 2n ,\] a contradiction. Without loss, suppose $G$ has a vertex, $x$ with in-degree one. Then \[ \rell{G}{p}=(1-p)\rell{G-x}{p}+p\rell{G|_{x}}{p},\] where $G|_{x}$ is the subdigraph of $G$ where $x$ must be operational.

First, if the degree of $x$ is 2 then $G-x$ has $2n-4=2(n-2)=2((n-1)-1)$ arcs and $n-1$ vertices so $\rell{S_{n-1}}{p}>\rell{G-x}{p}$ by our induction hypothesis. Otherwise, suppose for some value of $p\in (0,1)$ that $G-x$ is more optimal than $S_{n-1}$ and has fewer than $2n-4$ arcs. Then this implies that there exists a digraph on $2n-4$ arcs that is more optimal than $S_{n-1}$ for some value of $p\in (0,1)$, since we can add arcs to go from $G-x$ to some digraph $H$ on $2n-4$ arcs, such that it will still be strongly connected and every induced strongly connected subdigraph of $G-x$ will be strongly connected in $H$ as well. There may possibly exist even more such subdigraphs in $H$, which is a contradiction to our induction hypothesis. We now show that $\rell{G|_x}{p}\leq (1-p)^{n-1}+p$. 

Consider $G|_x$. Either $x$ is the only vertex operational, which occurs with probability $(1-p)^{n-1}$ or at least one other vertex is operational. Let \[\rell{G|_x}{p}=(1-p)^{n-1}+pP(G|_x),\] where $P(G|_x)$ is the probability that $x$ and the other operational vertices (of which there is at least one) induce a strongly connected subdigraph of $G$. We know that since $x$ has only one in-neighbour, call that vertex $y$, that must be operational as well. This occurs with probability $p$. There are two ways that $x$ and $y$ can be adjacent.

First, if $x$ and $y$ are adjacent via a bundle, then the optimal situation is that any remaining vertex can be operational or not, which means that $G$ is $S_{n-1}$ and $y$ is the central vertex and $P(G|_x)=1$, otherwise $P(G|_x)<1$. Secondly if $x$ and $y$ are not adjacent via bundle, then $x$ needs an operational out-neighbour to ensure that the resulting induced subdigraph is strongly connected, thus $P(G|_x)<1$.

This means
\begin{eqnarray*}
\rell{S_n}{p}-\rell{G}{p}&=&(1-p)[\rell{S_{n-1}}{p}-\rell{G}{p}]+\\
                 & & p[((1-p)^{n-1}+p)-(1-p)^{n-1}+pP(G_x)]\\
                 &=&(1-p)[\rell{S_{n-1}}{p}-\rell{G}{p}]+\\
                 & & p[p(1-P(G_x))]\\
                 &\geq &0\\
\end{eqnarray*}

Note that if $G$ is not $S_n$ then the inequality is strict, thus for $m=2n-2$, an optimally-greatest digraph exists and it is $S_n$.
\end{pf}

Consider the digraph, $D_n$ on $n$ vertices and $m=2n-1$ arcs which is $S_n$ with an additional arc between two of the degree two vertices. This digraph has the same reliability as $S_n$, as the additional arc is redundant, its removal does not affect the reliability of the digraph. For $n=3$ this is an optimally-greatest digraph. Following a very similar argument to Theorem~\ref{thrm:2n-2} we get the following result.

\begin{cor}
Let $G$ be a strongly connected digraph of order $n\geq 3$ and $m=2n-1$. Then $\rell{D_n}{p}\geq \rell{G}{p}$ for $p\in [0,1]$.
\end{cor}

\section{Digraphs with $n$ vertices and $2n$ arcs}

In this section we turn our focus to digraphs with $2n$ arcs. For node reliability an optimally-greatest graph on $n$ vertices and $n$ edges exists, namely the cycle. Of course, it is natural to think that the question of optimality for strongly connected node reliability may be trivial, just replace every edge of a cycle with a bundle, but we will show that there is no optimally-greatest digraph with $n$ vertices and $2n$ arcs. If such an optimally-greatest digraph did exist, then it would need to be optimal for values of $p$ sufficiently close to zero and for values of $p$ sufficiently close to one. Lemma \ref{Lemma: N compare} and Lemma \ref{Lemma: F compare} gives that for values of $p$ sufficiently close to 0, an optimally-greatest digraph will have the maximal number of bundles. For values of $p$ sufficiently close to one, a digraph with the largest vertex connectivity is most optimal.

A connected digraph $G$ with $n$ vertices and $2n$ arcs that has the most bundles is such that $U(G)$ is unicyclic as it will have $n$ edges. Of those digraphs, the one with the largest vertex connectivity is the cycle graph where every edge is replaced with a bundle. Therefore, to show that there is not an optimally-greatest digraph, we show that there is another digraph that is better than the bundled cycle graph for some $p\in (0,1)$. In particular, we show that the bundled cycle graph is not optimal for values of $p$ sufficiently close to one, there is another digraph that is more reliable.

Consider optimality for values of $p$ near 1 for a digraph $D$ with $n$ vertices and $2n$ arcs. Note that if one vertex has an in-degree (or out-degree) of one, then $D$ will not be most optimal for values of $p$ near 1 since it will have a cut vertex. Similarly, if $D$ has no vertices of in-degree (or out-degree) one, and has at least one vertex with in-degree (or out-degree) of at least three, then  $2n=\sum_{v\in V(D)}deg_+(v)\geq 2(n-1)+3=2n+1$, a contradiction to the number of arcs in the digraph. Therefore,  we can restrict our attention towards 2-regular digraphs.

As a bundled cycle belongs to the family of directed circulant graphs we turn our attention to studying the optimality of this family. Directed circulants are defined as follows. A \textit{generating set} of a group is a subset of the group set that does not contain the identity element such that every element of the group can be expressed as a combination under the group operation of finitely many elements of the subset and their inverses. Consider the cyclic group $(\mathbb{Z}_n,+_n)$ and let $S$ be a generating set of $(\mathbb{Z}_n,+_n)$.  We will label the vertices of a directed circulant graph of order $n$ as $v_0, v_1, \ldots ,v_{n-1}$. Then, the associated \textit{directed circulant graph}, denoted $\Gamma(\mathbb{Z}_n,S)$, has a vertex for each element of $G$, and there is a directed edge from $v_j$ to $v_i$, if and only if $i-j(\textrm{mod}~n)\in S$. Note, all arithmetic on the vertex subscripts is done modulo the order of the graph, and we omit the mod $n$ for ease of notation.

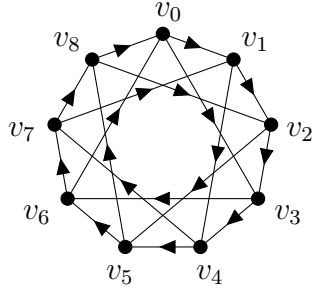
\begin{figure}[h!]
\centering
\begin{tikzpicture}[line cap=round,line join=round,>=triangle 45,x=1.0cm,y=1.0cm]
\draw [->-] (1.060307379214092,1.6275953626987478) -- (1.5603073792140922,2.493620766483186);
\draw [->-] (1.5603073792140922,2.493620766483186) -- (2.5,2.8356409098088537);
\draw [->-] (2.5,2.8356409098088537) -- (3.4396926207859084,2.4936207664831853);
\draw [->-] (3.4396926207859084,2.4936207664831853) -- (3.939692620785908,1.6275953626987465);
\draw [->-] (3.939692620785908,1.6275953626987465) -- (3.766044443118978,0.642787609686539);
\draw [->-] (3.766044443118978,0.642787609686539) -- (3.,0.);
\draw [->-] (3.,0.) -- (2.,0.);
\draw [->-] (2.,0.) -- (1.2339555568810219,0.6427876096865397);
\draw [->-] (1.2339555568810219,0.6427876096865397) -- (1.060307379214092,1.6275953626987478);
\draw [->-] (1.060307379214092,1.6275953626987478) -- (3.4396926207859084,2.4936207664831853);
\draw [->-] (3.4396926207859084,2.4936207664831853) -- (3.,0.);
\draw [->-] (3.,0.) -- (1.060307379214092,1.6275953626987478);
\draw [->-] (1.2339555568810219,0.6427876096865397) -- (2.5,2.8356409098088537);
\draw [->-] (2.5,2.8356409098088537) -- (3.766044443118978,0.6427876096865388);
\draw [->-] (3.766044443118978,0.642787609686539) -- (1.233955556881022,0.6427876096865397);
\draw [->-] (2.,0.) -- (1.5603073792140922,2.493620766483186);
\draw [->-] (1.5603073792140922,2.493620766483186) -- (3.939692620785908,1.6275953626987465);
\draw [->-] (3.939692620785908,1.6275953626987465) -- (2.,0.);
\draw (1.6,-.1) node[anchor=north west] {$v_5$};
\draw (0.5,0.7) node[anchor=north west] {$v_6$};
\draw (0.3,1.8) node[anchor=north west] {$v_7$};
\draw (.9,3) node[anchor=north west] {$v_8$};
\draw (2.25,3.4) node[anchor=north west] {$v_0$};
\draw (3.4,3) node[anchor=north west] {$v_1$};
\draw (4,1.8) node[anchor=north west] {$v_2$};
\draw (3.84,0.7) node[anchor=north west] {$v_3$};
\draw (2.8,-.1) node[anchor=north west] {$v_4$};
\begin{scriptsize}
\draw [fill=black] (2.,0.) circle (2.5pt);
\draw [fill=black] (3.,0.) circle (2.5pt);
\draw [fill=black] (3.766044443118978,0.642787609686539) circle (2.5pt);
\draw [fill=black] (3.939692620785908,1.6275953626987465) circle (2.5pt);
\draw [fill=black] (3.4396926207859084,2.4936207664831853) circle (2.5pt);
\draw [fill=black] (2.5,2.8356409098088537) circle (2.5pt);
\draw [fill=black] (1.5603073792140922,2.493620766483186) circle (2.5pt);
\draw [fill=black] (1.060307379214092,1.6275953626987478) circle (2.5pt);
\draw [fill=black] (1.2339555568810219,0.6427876096865397) circle (2.5pt);
\end{scriptsize}
\end{tikzpicture}
\caption{Digraph $\Gamma(\mathbb{Z}_9,\{1,3\})$}
\label{Figure: Example 1}
\end{figure}

As an example, consider the additive group $\mathbb{Z}_9$ with group operation $+_9$. Let a generating set of $(\mathbb{Z}_9,+_9)$ be the set $S=\{1,3\}$. Then, there is a directed edge from $v_j$ to $v_i$ if and only if $(i-j)(\textrm{mod}~9)=1$ or $(i-j)(\textrm{mod}~9)=3$. The circulant digraph $\Gamma(\mathbb{Z}_9,\{1,3\})$ illustrating these adjacencies is shown in Figure \ref{Figure: Example 1}. 

By studying directed circulant graphs, we will show that an optimally-greatest digraph with $n$ vertices and $2n$ arcs does not exist. Furthermore, we will prove which directed circulant graph is optimal for values of $p$ near one.

There are three properties of circulant graphs that we use throughout this paper. The first gives a condition to determine which circulant graphs are isomorphic. This is done by using an \textit{automorphism of a group}, which is an isomorphism from a group to itself. It is well known that the automorphisms of $\mathbb{Z}_n$ are of the form $\phi_a:\mathbb{Z}_n\to\mathbb{Z}_n$, where $\phi_a(1)=a$ and $a\in\mathbb{Z}_n^*$, the multiplicative group of $\mathbb{Z}_n$. 

\begin{lem}[\cite{muzychuk1997adam}]\label{Lemma: Circ Iso}
If $S$ and $S'$ are related by a group automorphism, then $\Gamma(\mathbb{Z}_n,S)$ and $\Gamma(\mathbb{Z}_n,S')$ are isomorphic.
\end{lem}

The second property determines which circulant graphs are disconnected. This is done by checking the greatest common divisor of the elements in the generating set as well as the order of the graph.

\begin{lem}[\cite{dummit1991abstract}]\label{Lemma: Circ Disconnected}
If $\gcd(n,a,b)>1$, then $\Gamma(\mathbb{Z}_n,\{a,b\})$ is disconnected.
\end{lem}

The last property gives that the two in-neighbours of a vertex are also the two out-neighbours of another vertex. This result follows directly from the definition of the circulant graph $\Gamma(\mathbb{Z}_n,\{a,b\})$.

\begin{lem}\label{Lemma: Circulant Out In} 
For the digraph $\Gamma(\mathbb{Z}_n,\{a,b\})$, the vertices $v_{i+a}$ and $v_{i+b}$ are both the out-neighbours of the vertex $v_{i}$ and are both the in-neighbours of the vertex $v_{i+a+b}$.
\end{lem}

We will now determine which circulant graph of the form $\Gamma\left(\mathbb{Z}_n,\{a,b\}\right)$ is optimally-greatest for values of $p$ near one. The form of the circulant graph is dependent on the following conditions: If the order of the graph is even, if the order is odd and not divisible by three, and lastly if the order is odd and is divisible by three. 

\subsection{The order of the digraph is even}

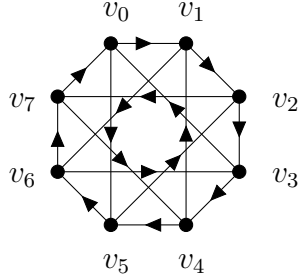
\begin{figure}[h!]
\centering
\begin{tikzpicture}[line cap=round,line join=round,>=triangle 45,x=1.0cm,y=1.0cm]
\draw [->-] (1.,2.414213562373095)-- (2.,2.414213562373095);
\draw [->-] (2.7071067811865475,1.7071067811865472)-- (2.7071067811865475,0.7071067811865475);
\draw [->-] (2.7071067811865475,0.7071067811865475)-- (2.,0.);
\draw [->-] (2.,0.)-- (1.,0.);
\draw [->-] (1.,0.)-- (0.2928932188134523,0.7071067811865478);
\draw [->-] (0.2928932188134523,0.7071067811865478)-- (0.29289321881345254,1.7071067811865477);
\draw [->-] (0.29289321881345254,1.7071067811865477)-- (1.,2.414213562373095);
\draw [->-] (2.,2.414213562373095)-- (2.7071067811865475,1.7071067811865472);
\draw [->-] (1.,2.414213562373095)-- (1.,0.);
\draw [->-] (2.,2.414213562373095)-- (0.2928932188134523,0.7071067811865478);
\draw [->-] (2.7071067811865475,1.7071067811865472)-- (0.29289321881345254,1.7071067811865477);
\draw [->-] (2.7071067811865475,0.7071067811865475)-- (1.,2.414213562373095);
\draw [->-] (2.,0.)-- (2.,2.414213562373095);
\draw [->-] (1.,0.)-- (2.7071067811865475,1.7071067811865472);
\draw [->-] (0.2928932188134523,0.7071067811865478)-- (2.7071067811865475,0.7071067811865475);
\draw [->-] (0.29289321881345254,1.7071067811865477)-- (2.,0.);
\draw (0.75,3.1) node[anchor=north west] {$v_0$};
\draw (1.75,3.1) node[anchor=north west] {$v_1$};
\draw (3,1.9) node[anchor=north west] {$v_2$};
\draw (3,0.9) node[anchor=north west] {$v_3$};
\draw (1.75,-0.25) node[anchor=north west] {$v_4$};
\draw (0.75,-0.25) node[anchor=north west] {$v_5$};
\draw (-0.5,0.9) node[anchor=north west] {$v_6$};
\draw (-0.5,1.9) node[anchor=north west] {$v_7$};
\begin{scriptsize}
\draw [fill=black] (1.,0.) circle (2.5pt);
\draw [fill=black] (2.,0.) circle (2.5pt);
\draw [fill=black] (2.7071067811865475,0.7071067811865475) circle (2.5pt);
\draw [fill=black] (2.7071067811865475,1.7071067811865472) circle (2.5pt);
\draw [fill=black] (2.,2.414213562373095) circle (2.5pt);
\draw [fill=black] (1.,2.414213562373095) circle (2.5pt);
\draw [fill=black] (0.29289321881345254,1.7071067811865477) circle (2.5pt);
\draw [fill=black] (0.2928932188134523,0.7071067811865478) circle (2.5pt);
\end{scriptsize}
\end{tikzpicture}
\caption{Digraph $\Gamma(\mathbb{Z}_8,\{1,5\})$}
\label{Figure: Circ8}
\end{figure}

In this section, we will show that if the order of the digraph is even, then the digraph $\Gamma\left(\mathbb{Z}_{2k},\{1,k+1\}\right)$ is optimally-greatest for values of $p$ sufficiently close to one. For illustration, the digraph $\Gamma\left(\mathbb{Z}_{8},\{1,5\}\right)$ is provided in Figure \ref{Figure: Circ8}

Notice for any digraph, if an operational vertex has no operational out-neighbours or no operational in-neighbours, then the resulting induced subdigraph is not strongly connected. An induced subdigraph formed this way is called a \textit{trivial subdigraph}. Then, since Lemma \ref{Lemma: Circulant Out In} gives that the in-neighbours of a vertex are the out-neighbours of another vertex, in counting the number of trivial subdigraphs, it is sufficient to either consider just the out-neighbours of the vertices or just the in-neighbours of the vertices. In this paper, we consider the out-neighbours of the vertices. 

The following lemma shows that for $\Gamma\left(\mathbb{Z}_{2k},\{1,k+1\}\right)$ to have an induced subdigraph that is not strongly connected when two vertices are not operational the two non-operational vertices must both the out-neighbours of a vertex. That is, the subdigraphs of $\Gamma\left(\mathbb{Z}_{2k},\{1,k+1\}\right)$, in this case, that are not strongly connected are the trivial subdigraphs.  

\begin{lem}
\label{Lemma: Property Even}
If two or fewer vertices are not operational and all of the  operational vertices have at least one operational out-neighbour and one operational in-neighbour, then the resulting induced subdigraph of $\Gamma\left(\mathbb{Z}_{2k},\{1,k+1\}\right)$ is strongly connected.
\end{lem}
\begin{pf}
Without loss of generality, suppose that there are two vertices that are not operational and that $v_0$ is an operational vertex. Let $i\neq 0$ be the smallest index such that $v_i$ is not operational. Then, there are two cases to consider. The first case is where $v_{i+1}$ is not operational and the second case is where $v_{i+1}$ is operational.

Firstly, suppose that $v_{i+1}$ is not operational. Then, since all of the  operational vertices have at least one operational out-neighbour, the vertices $v_{i+k}$ and $v_{i+k+1}$ are operational. This is because $v_{i}$ and $v_{i+k}$ are both out-neighbours of $v_{i-1}$ and because $v_{i+1}$ and $v_{i+k+1}$ are both out-neighbours of $v_{i+k}$ (Lemma~\ref{Lemma: Circulant Out In}). Now, consider the following sequence of vertices.
$$\{v_0,v_1,\dots,v_{i-1},v_{i+k},v_{i+k+1},v_{i+2},v_{i+3},\dots,v_{2k-1},v_0\}$$

This sequence of vertices provides a directed circuit containing all the operational vertices of the induced subdigraph. Therefore, when $v_i$ and $v_{i+1}$ are not operational, the induced subdigraph is strongly connected. For illustration, the case where $v_2$ and $v_3$ are not operational for $\Gamma(\mathbb{Z}_8,\{1,5\})$ is given in Figure \ref{Figure: Circ8 Case 1}. The directed path from $v_0$ to $v_4$ is highlighted in red.

\begin{figure}[h!]
\centering
\begin{tikzpicture}[line cap=round,line join=round,>=triangle 45,x=1.0cm,y=1.0cm]
\draw [->-,color=red] (1.,2.414213562373095)-- (2.,2.414213562373095);
\draw [->-] (2.,0.)-- (1.,0.);
\draw [->-] (1.,0.)-- (0.2928932188134523,0.7071067811865478);
\draw [->-,color=red] (0.2928932188134523,0.7071067811865478)-- (0.29289321881345254,1.7071067811865477);
\draw [->-] (0.29289321881345254,1.7071067811865477)-- (1.,2.414213562373095);
\draw [->-] (1.,2.414213562373095)-- (1.,0.);
\draw [->-,color=red] (2.,2.414213562373095)-- (0.2928932188134523,0.7071067811865478);
\draw [->-] (2.,0.)-- (2.,2.414213562373095);
\draw [->-,color=red] (0.29289321881345254,1.7071067811865477)-- (2.,0.);
\draw (0.75,3.1) node[anchor=north west] {$v_0$};
\draw (1.75,3.1) node[anchor=north west] {$v_1$};
\draw (1.75,-0.25) node[anchor=north west] {$v_4$};
\draw (0.75,-0.25) node[anchor=north west] {$v_5$};
\draw (-0.5,0.9) node[anchor=north west] {$v_6$};
\draw (-0.5,1.9) node[anchor=north west] {$v_7$};
\begin{scriptsize}
\draw [fill=black] (1.,0.) circle (2.5pt);
\draw [fill=black] (2.,0.) circle (2.5pt);
\draw [fill=black] (2.,2.414213562373095) circle (2.5pt);
\draw [fill=black] (1.,2.414213562373095) circle (2.5pt);
\draw [fill=black] (0.29289321881345254,1.7071067811865477) circle (2.5pt);
\draw [fill=black] (0.2928932188134523,0.7071067811865478) circle (2.5pt);
\end{scriptsize}
\end{tikzpicture}
\caption{Subdigraph 1 of $\Gamma(\mathbb{Z}_8,\{1,5\})$}
\label{Figure: Circ8 Case 1}
\end{figure}

Now, suppose $v_{i+1}$ is operational. Then, there is some vertex $v_j$, where $j>i+1$, such that $v_j$ is not operational. Since all of the operational vertices have at least one operational out-neighbour, the vertices $v_{i+k}$ and $v_{j+k}$ are operational. This is because $v_i$ and $v_{i+k}$ are both out-neighbours of $v_{i-1}$ and $v_{j+k}$ and $v_j$ are both out-neighbours of $v_{j-1}$. Now, consider the following sequence of vertices.
$$\{v_0,v_1,\dots,v_{i-1},v_{i+k},v_{i+1},v_{i+2},\dots,v_{j-2},v_{j-1},v_{j+k},v_{j+1},v_{j+2},\dots, v_{2k-1},v_0\}$$

This sequence of vertices provides a directed circuit containing all the operational vertices of the induced subdigraph. Therefore, when $v_i$ and $v_j$, $j\neq i+1$, are not operational, the induced subdigraph is strongly connected. For illustration, the case where $v_1$ and $v_3$ are not operational for $\Gamma(\mathbb{Z}_8,\{1,5\})$ is given in Figure \ref{Figure: Circ8 Case 2}. The directed path from $v_0$ to $v_4$ is highlighted in red.

\begin{figure}[h!]
\centering
\begin{tikzpicture}[line cap=round,line join=round,>=triangle 45,x=1.0cm,y=1.0cm]
\draw [->-] (2.,0.)-- (1.,0.);
\draw [->-] (1.,0.)-- (0.2928932188134523,0.7071067811865478);
\draw [->-] (0.2928932188134523,0.7071067811865478)-- (0.29289321881345254,1.7071067811865477);
\draw [->-] (0.29289321881345254,1.7071067811865477)-- (1.,2.414213562373095);
\draw [->-,color=red] (1.,2.414213562373095)-- (1.,0.);
\draw [->-,color=red] (2.7071067811865475,1.7071067811865472)-- (0.29289321881345254,1.7071067811865477);
\draw [->-,color=red] (1.,0.)-- (2.7071067811865475,1.7071067811865472);
\draw [->-,color=red] (0.29289321881345254,1.7071067811865477)-- (2.,0.);
\draw (0.75,3.1) node[anchor=north west] {$v_0$};
\draw (3,1.9) node[anchor=north west] {$v_2$};
\draw (1.75,-0.25) node[anchor=north west] {$v_4$};
\draw (0.75,-0.25) node[anchor=north west] {$v_5$};
\draw (-0.5,0.9) node[anchor=north west] {$v_6$};
\draw (-0.5,1.9) node[anchor=north west] {$v_7$};
\begin{scriptsize}
\draw [fill=black] (1.,0.) circle (2.5pt);
\draw [fill=black] (2.,0.) circle (2.5pt);
\draw [fill=black] (2.7071067811865475,1.7071067811865472) circle (2.5pt);
\draw [fill=black] (1.,2.414213562373095) circle (2.5pt);
\draw [fill=black] (0.29289321881345254,1.7071067811865477) circle (2.5pt);
\draw [fill=black] (0.2928932188134523,0.7071067811865478) circle (2.5pt);
\end{scriptsize}
\end{tikzpicture}
\caption{Subdigraph 2 of $\Gamma(\mathbb{Z}_8,\{1,5\})$}
\label{Figure: Circ8 Case 2}
\end{figure}

Therefore, in both cases, the resulting induced subdigraph of $\Gamma\left(\mathbb{Z}_{2k},\{1,k+1\}\right)$ is strongly connected. As the connectivity of  $\Gamma\left(\mathbb{Z}_{2k},\{1,k+1\}\right)$ is 2, the only induced subdigraphs that are not strongly connected are the trivial subdigraphs.
\end{pf}

The proof of the following optimality theorem is as follows. First, Lemma \ref{Lemma: Property Even} will be applied to give the exact $F_1$ and $F_2$ values for $\Gamma\left(\mathbb{Z}_{2k},\{1,k+1\}\right)$. Then, we will show that while some other circulant digraphs may have the same $F_1$ term, the $F_2$ term for $\Gamma\left(\mathbb{Z}_{2k},\{1,k+1\}\right)$ is larger. 

\begin{thrm}
\label{Theorem: Even Circulant}
If $G=\Gamma\left(\mathbb{Z}_{2k},\{1,k+1\}\right)$, then $G$ is an optimally-greatest circulant digraph of the form $\Gamma\left(\mathbb{Z}_{2k},\{a,b\}\right)$ for values of $p$ near one.
\end{thrm}
\begin{pf}
Suppose first that only one vertex of $G$ is not operational. By Lemma \ref{Lemma: Property Even}, the resulting induced subdigraph of $G$ is strongly connected. Therefore, the $F_1$ term of $G$ is maximal. Now, suppose that two vertices of $G$ are not operational. By Lemma \ref{Lemma: Property Even}, if all of the  operational vertices of $G$ have at least one operational in-neighbour and one operational out-neighbour, then the resulting induced subdigraph of $G$ is strongly connected. Thus by Lemma \ref{Lemma: Circulant Out In}, the remaining cases left to consider are how many ways there are to remove two out-neighbours of a vertex.

The two out-neighbours of a vertex $v_i$ are $v_{i+1}$ and $v_{i+k+1}$. However, notice that $v_{i+1}$ and $v_{i+k+1}$ are also out-neighbours of $v_{i+k}$. This is because $(i+k+1)-(i+k)=1$ and $(i+1)-(i+k)=1-k\equiv 1+k$(mod $2k$). Thus, there are only $k$ ways to remove two out-neighbours of a vertex. This is illustrated for $\Gamma(\mathbb{Z}_8,\{1,5\})$ in Figure \ref{Figure: Circ8 Double Arcs}.

\begin{figure}[h!]
\centering
\begin{tikzpicture}[line cap=round,line join=round,>=triangle 45,x=1.0cm,y=1.0cm]
\draw [->-,color=red] (1.,2.414213562373095)-- (2.,2.414213562373095);
\draw [->-] (2.7071067811865475,1.7071067811865472)-- (2.7071067811865475,0.7071067811865475);
\draw [->-] (2.7071067811865475,0.7071067811865475)-- (2.,0.);
\draw [->-,color=red] (2.,0.)-- (1.,0.);
\draw [->-] (1.,0.)-- (0.2928932188134523,0.7071067811865478);
\draw [->-] (0.2928932188134523,0.7071067811865478)-- (0.29289321881345254,1.7071067811865477);
\draw [->-] (0.29289321881345254,1.7071067811865477)-- (1.,2.414213562373095);
\draw [->-] (2.,2.414213562373095)-- (2.7071067811865475,1.7071067811865472);
\draw [->-,color=red] (1.,2.414213562373095)-- (1.,0.);
\draw [->-] (2.,2.414213562373095)-- (0.2928932188134523,0.7071067811865478);
\draw [->-] (2.7071067811865475,1.7071067811865472)-- (0.29289321881345254,1.7071067811865477);
\draw [->-] (2.7071067811865475,0.7071067811865475)-- (1.,2.414213562373095);
\draw [->-,color=red] (2.,0.)-- (2.,2.414213562373095);
\draw [->-] (1.,0.)-- (2.7071067811865475,1.7071067811865472);
\draw [->-] (0.2928932188134523,0.7071067811865478)-- (2.7071067811865475,0.7071067811865475);
\draw [->-] (0.29289321881345254,1.7071067811865477)-- (2.,0.);
\draw (0.75,3.1) node[anchor=north west] {$v_0$};
\draw (1.75,3.1) node[anchor=north west] {$v_1$};
\draw (3,1.9) node[anchor=north west] {$v_2$};
\draw (3,0.9) node[anchor=north west] {$v_3$};
\draw (1.75,-0.25) node[anchor=north west] {$v_4$};
\draw (0.75,-0.25) node[anchor=north west] {$v_5$};
\draw (-0.5,0.9) node[anchor=north west] {$v_6$};
\draw (-0.5,1.9) node[anchor=north west] {$v_7$};
\begin{scriptsize}
\draw [fill=black] (1.,0.) circle (2.5pt);
\draw [fill=black] (2.,0.) circle (2.5pt);
\draw [fill=black] (2.7071067811865475,0.7071067811865475) circle (2.5pt);
\draw [fill=black] (2.7071067811865475,1.7071067811865472) circle (2.5pt);
\draw [fill=black] (2.,2.414213562373095) circle (2.5pt);
\draw [fill=black] (1.,2.414213562373095) circle (2.5pt);
\draw [fill=black] (0.29289321881345254,1.7071067811865477) circle (2.5pt);
\draw [fill=black] (0.2928932188134523,0.7071067811865478) circle (2.5pt);
\end{scriptsize}
\end{tikzpicture}
\caption{$\Gamma(\mathbb{Z}_8,\{1,5\})$}
\label{Figure: Circ8 Double Arcs}
\end{figure}

Since there are $k(2k-1)$ ways to have two non-operational vertices, and only $k$ of these result in an induced subdigraph that is not strongly connected, the $F_2$ term for $G$ is given by $F_2(G)=k(2k-1)-k=k(2k-2)$. Lastly, it remains to show that no other non-isomorphic directed circulant graphs have the property that the two out-neighbours of a vertex are also out-neighbours of another vertex. 

Consider a circulant digraph $\Gamma(\mathbb{Z}_{2k},\{a,b\})$ and the vertices $v_a$ and $v_b$. These two vertices are both out-neighbours of the vertex $v_0$. Now, suppose $v_a$ and $v_b$ are also out-neighbours of the vertex $v_c$. Then, $a-c\equiv b$(mod $2k$) and $b-c \equiv a$(mod $2k$), which gives that $a-b\equiv b-a$(mod $2k$), which occurs when $a-b\equiv k$(mod $2k$). Hence, the directed circulant graphs to consider are of the form $\Gamma(\mathbb{Z}_{2k},\{c,c+k\})$. 

We will now show that all of the non-isomorphic circulant digraphs of the form $\Gamma(\mathbb{Z}_{2k},\{c,c+k\})$ are disconnected. First, suppose that $c$ is a unit, which implies that $c$ must be odd. Then, multiplying the elements of the set $\{1,k+1\}$ by $c$ gives the resulting set $\{c,ck+c\}\equiv \{c,k+c\}$(mod $2k$). Therefore, when $c$ is a unit, $\Gamma(\mathbb{Z}_{2k},\{c,c+k\})$ and $\Gamma\left(\mathbb{Z}_{2k},\{1,k+1\}\right)$ are isomorphic by Lemma \ref{Lemma: Circ Iso}. 

Next, suppose instead that $c+k$ is a unit, which implies that $c+k+1$ is even. Then, multiplying the elements of the set $\{1,k+1\}$ by $c+k$ gives the resulting set $\{c+k,(c+k)(k+1)\}\equiv \{c+k,c+(c+k+1)k\}\equiv \{c+k,c\}$(mod $2k$). Therefore, when $c+k$ is a unit, $\Gamma(\mathbb{Z}_{2k},\{c,c+k \})$ and $\Gamma\left(\mathbb{Z}_{2k},\{1,k+1\}\right)$ are isomorphic by Lemma \ref{Lemma: Circ Iso}.

Recall that Lemma \ref{Lemma: Circ Disconnected} states that if $\gcd(n,a,b)\neq 1$, then the circulant digraph $\Gamma(\mathbb{Z}_{n},\{a,b\})$ is disconnected. Therefore, we will show that when both $c$ and $c+k$ are not units, that is when $c$ shares a factor with $2k$ and $c+k$ shares a factor with $2k$, that this condition holds. Thus, the directed circulant graphs of this form are disconnected. 

If $c$ is odd, then from properties of the greatest common divisor, it follows that $1<\gcd(2k,c)=\gcd(k,c)$. Thus, it must be the case that $c$ and $k$ share a common factor $r>1$. Therefore, $\gcd(2k,c,c+k)\geq r>1$, which gives that the circulant graph is disconnected.

If $c$ is even and $k$ is even, then $\gcd(2k,c,c+k)\geq 2$, which gives that the circulant graph is disconnected. If $c$ is even and $k$ is odd, then from properties of the greatest common divisor, it follows that $1<\gcd(c+k,2k)=\gcd(c+k,k)=\gcd(c,k)$. Thus, $c$ and $k$ share a common factor $r>1$. Therefore, $\gcd(2k,c,c+k)\geq r>1$, which gives that the circulant graph is disconnected. 
\end{pf}

\subsection{The order of the graph is odd and not divisible by three}

\begin{figure}[h!]
\centering
\begin{tikzpicture}[line cap=round,line join=round,>=triangle 45,x=1.0cm,y=1.0cm]
\draw [->-] (1.04,-0.02)--(0.,0.) ;
\draw [->-] (1.9257164904935398,0.5254413794971975)--(1.04,-0.02) ;
\draw [->-] (2.375940743922592,1.4631503544058582)--(1.9257164904935398,0.5254413794971975) ;
\draw [->-](2.2477297409559944,2.4954109507274937)-- (2.375940743922592,1.4631503544058582);
\draw [->-] (1.5817895691399835,3.294487722734828)--(2.2477297409559944,2.4954109507274937) ;
\draw [->-] (0.5895515277177357,3.606679441322205)--(1.5817895691399835,3.294487722734828) ;
\draw [->-] (-0.4139558159781702,3.332867441679409)-- (0.5895515277177357,3.606679441322205);
\draw [->-] (-1.1101259707683522,2.559985099029886)--(-0.4139558159781702,3.332867441679409) ;
\draw [->-] (-1.2779298314101877,1.5334170962391824)--(-1.1101259707683522,2.559985099029886) ;
\draw [->-] (-0.8640908577953159,0.579091520810445)--(-1.2779298314101877,1.5334170962391824) ;
\draw [->-] (0.,0.)--(-0.8640908577953159,0.579091520810445) ;
\draw [->-] (0.5895515277177357,3.606679441322205)-- (-1.2779298314101877,1.5334170962391824);
\draw [->-] (1.5817895691399835,3.294487722734828)-- (-1.1101259707683522,2.559985099029886);
\draw [->-] (2.2477297409559944,2.4954109507274937)-- (-0.4139558159781702,3.332867441679409);
\draw [->-] (2.375940743922592,1.4631503544058582)-- (0.5895515277177357,3.606679441322205);
\draw [->-] (1.9257164904935398,0.5254413794971975)-- (1.5817895691399835,3.294487722734828);
\draw [->-] (1.04,-0.02)-- (2.2477297409559944,2.4954109507274937);
\draw [->-] (0.,0.)-- (2.375940743922592,1.4631503544058582);
\draw [->-] (-0.8640908577953159,0.579091520810445)-- (1.9257164904935398,0.5254413794971975);
\draw [->-] (-1.2779298314101877,1.5334170962391824)-- (1.04,-0.02);
\draw [->-] (-1.1101259707683522,2.559985099029886)-- (0.,0.);
\draw [->-] (-0.4139558159781702,3.332867441679409)-- (-0.8640908577953159,0.579091520810445);
\draw (0.35,4.22) node[anchor=north west] {$v_0$};
\draw (1.75,3.7) node[anchor=north west] {$v_1$};
\draw (2.4,2.75) node[anchor=north west] {$v_2$};
\draw (2.45,1.65) node[anchor=north west] {$v_3$};
\draw (1.85,0.5) node[anchor=north west] {$v_4$};
\draw (.85,-0.15) node[anchor=north west] {$v_5$};
\draw (-.35,-0.15) node[anchor=north west] {$v_6$};
\draw (-1.4,0.5) node[anchor=north west] {$v_7$};
\draw (-2,1.65) node[anchor=north west] {$v_8$};
\draw (-1.85,2.75) node[anchor=north west] {$v_9$};
\draw (-1.4,3.7) node[anchor=north west] {$v_{10}$};
\begin{scriptsize}
\draw [fill=black] (1.04,-0.02) circle (2.5pt);
\draw [fill=black] (1.9257164904935398,0.5254413794971975) circle (2.0pt);
\draw [fill=black] (2.375940743922592,1.4631503544058582) circle (2.0pt);
\draw [fill=black] (2.2477297409559944,2.4954109507274937) circle (2.0pt);
\draw [fill=black] (1.5817895691399835,3.294487722734828) circle (2.0pt);
\draw [fill=black] (0.5895515277177357,3.606679441322205) circle (2.0pt);
\draw [fill=black] (-0.4139558159781702,3.332867441679409) circle (2.0pt);
\draw [fill=black] (-1.1101259707683522,2.559985099029886) circle (2.0pt);
\draw [fill=black] (-1.2779298314101877,1.5334170962391824) circle (2.0pt);
\draw [fill=black] (-0.8640908577953159,0.579091520810445) circle (2.0pt);
\draw [fill=black] (0.,0.) circle (2.0pt);
\end{scriptsize}
\end{tikzpicture}
\caption{Digraph $\Gamma(\mathbb{Z}_{11},\{1,8\})$}
\label{Figure: Circ11}
\end{figure}
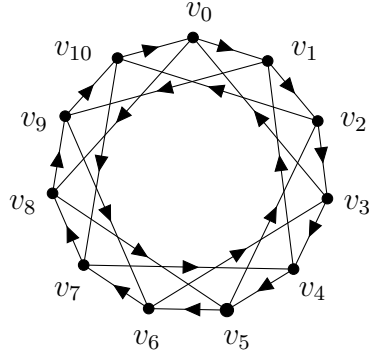

In this section, we will show that if $n$ is odd and not divisible by three, then $\Gamma(\mathbb{Z}_n,\{1,2(3^{-1}))$ is an optimally-greatest circulant digraph for values of $p$ sufficiently close to one. For illustration, the graph $\Gamma(\mathbb{Z}_{11},\{1,8\})$ is provided in Figure \ref{Figure: Circ11}. 

We show that if five or fewer vertices are not operational, having an operational vertex with no operational out-neighbours is the only way for $\Gamma\left(\mathbb{Z}_{n},\{1,2(3^{-1})\}\right)$ to have an induced subdigraph that is not strongly connected. That is, the only induced subdigraphs that are not strongly connected in this case are the trivial subdigraphs.

\begin{lem}
\label{Lemma: Property Odd Not Divisible By Three}
If five or fewer vertices are not operational and all of the  operational vertices have at least one operational out-neighbour and one operational in-neighbour, then the resulting induced subdigraph of $\Gamma\left(\mathbb{Z}_{n},\{1,2(3^{-1})\}\right)$ is strongly connected.
\end{lem}
\begin{pf}
Without loss of generality, suppose that five vertices are not operational and that $v_0$ is operational. Let $i\neq 0$ be the smallest index such that $v_i$ is not operational. Then, since all of the  operational vertices have at least one operational out-neighbour, the vertex $v_{i-1+2(3^{-1})}$ is operational. This is because the vertices $v_{i}$ and $v_{i-1+2(3^{-1})}$ are both out-neighbours of $v_{i-1}$. 

\textbf{Case 1:} First, suppose that the vertex $v_{i-1+4(3^{-1})}$ is operational. Therefore, the following directed path exists.
$$\{v_0,v_1,\dots,v_{i-1},v_{i-1+2(3^{-1})},v_{i-1+4(3^{-1})}\}$$

Now, there are two cases to consider. Either the vertex $v_{i+1}$ is operational or it is not operational. Suppose first that $v_{i+1}$ is operational. Then, the following directed path exists.$$\{v_0,v_1,\dots,v_{i-1},v_{i-1+2(3^{-1})},v_{i-1+4(3^{-1})},v_{i+1}\}$$ 

For illustration, the case where $v_2$ is not operational for the directed circulant graph $\Gamma(\mathbb{Z}_{11},\{1,8\})$ is given in Figure \ref{Figure: Circ11 Case1}. Also, the directed path from the vertex $v_0$ to the vertex $v_3$ is highlighted in red.

\begin{figure}[h!]
\centering
\begin{tikzpicture}[line cap=round,line join=round,>=triangle 45,x=1.0cm,y=1.0cm]
\draw [->-] (1.04,-0.02)--(0.,0.) ;
\draw [->-] (1.9257164904935398,0.5254413794971975)--(1.04,-0.02) ;
\draw [->-] (2.375940743922592,1.4631503544058582)--(1.9257164904935398,0.5254413794971975) ;
\draw [->-,color=red] (0.5895515277177357,3.606679441322205)--(1.5817895691399835,3.294487722734828) ;
\draw [->-] (-0.4139558159781702,3.332867441679409)-- (0.5895515277177357,3.606679441322205);
\draw [->-] (-1.1101259707683522,2.559985099029886)--(-0.4139558159781702,3.332867441679409) ;
\draw [->-] (-1.2779298314101877,1.5334170962391824)--(-1.1101259707683522,2.559985099029886) ;
\draw [->-] (-0.8640908577953159,0.579091520810445)--(-1.2779298314101877,1.5334170962391824) ;
\draw [->-] (0.,0.)--(-0.8640908577953159,0.579091520810445) ;

\draw [->-] (0.5895515277177357,3.606679441322205)-- (-1.2779298314101877,1.5334170962391824);
\draw [->-,color=red] (1.5817895691399835,3.294487722734828)-- (-1.1101259707683522,2.559985099029886);
\draw [->-] (2.375940743922592,1.4631503544058582)-- (0.5895515277177357,3.606679441322205);
\draw [->-] (1.9257164904935398,0.5254413794971975)-- (1.5817895691399835,3.294487722734828);
\draw [->-,color=red] (0.,0.)-- (2.375940743922592,1.4631503544058582);
\draw [->-] (-0.8640908577953159,0.579091520810445)-- (1.9257164904935398,0.5254413794971975);
\draw [->-] (-1.2779298314101877,1.5334170962391824)-- (1.04,-0.02);
\draw [->-,color=red] (-1.1101259707683522,2.559985099029886)-- (0.,0.);
\draw [->-] (-0.4139558159781702,3.332867441679409)-- (-0.8640908577953159,0.579091520810445);
\draw (0.35,4.22) node[anchor=north west] {$v_0$};
\draw (1.75,3.7) node[anchor=north west] {$v_1$};
\draw (2.45,1.65) node[anchor=north west] {$v_3$};
\draw (1.85,0.5) node[anchor=north west] {$v_4$};
\draw (.85,-0.15) node[anchor=north west] {$v_5$};
\draw (-.35,-0.15) node[anchor=north west] {$v_6$};
\draw (-1.4,0.5) node[anchor=north west] {$v_7$};
\draw (-2,1.65) node[anchor=north west] {$v_8$};
\draw (-1.85,2.75) node[anchor=north west] {$v_9$};
\draw (-1.4,3.7) node[anchor=north west] {$v_{10}$};
\begin{scriptsize}
\draw [fill=black] (1.04,-0.02) circle (2.5pt);
\draw [fill=black] (1.9257164904935398,0.5254413794971975) circle (2.0pt);
\draw [fill=black] (2.375940743922592,1.4631503544058582) circle (2.0pt);
\draw [fill=black] (1.5817895691399835,3.294487722734828) circle (2.0pt);
\draw [fill=black] (0.5895515277177357,3.606679441322205) circle (2.0pt);
\draw [fill=black] (-0.4139558159781702,3.332867441679409) circle (2.0pt);
\draw [fill=black] (-1.1101259707683522,2.559985099029886) circle (2.0pt);
\draw [fill=black] (-1.2779298314101877,1.5334170962391824) circle (2.0pt);
\draw [fill=black] (-0.8640908577953159,0.579091520810445) circle (2.0pt);
\draw [fill=black] (0.,0.) circle (2.0pt);
\end{scriptsize}
\end{tikzpicture}
\caption{Subdigraph 1 of $\Gamma(\mathbb{Z}_{11},\{1,8\})$}
\label{Figure: Circ11 Case1}
\end{figure}

Next, suppose now that $v_{i+1}$ is not operational and that $v_{i+2}$ is operational. Then, since each operational vertex has an operational out-neighbour, the vertex $v_{i+4(3^{-1})}$ is operational. Therefore, the following directed path exists $$\{v_0,v_1,\dots,v_{i-1},v_{i-1+2(3^{-1})},v_{i-1+4(3^{-1})},v_{i+4(3^{-1})},v_{i+2}\}$$ For illustration, the case where $v_2$ and $v_3$ are not operational for $\Gamma(\mathbb{Z}_{11},\{1,8\})$ is given in Figure \ref{Figure: Circ11 Case2}. Also, the directed path from $v_0$ to $v_4$ is highlighted in red. Note that the same arguments  could be applied reiteratively if $v_i,v_{i+1},v_{i+2},v_{i+3}$, and $v_{i+4}$ were all non-operational vertices. 

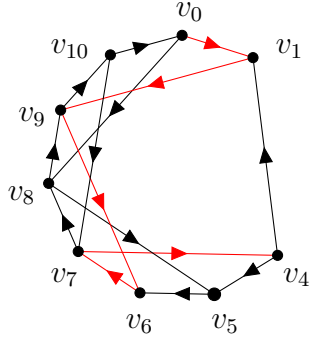
\begin{figure}[h!]
\centering
\begin{tikzpicture}[scale=.95,line cap=round,line join=round,>=triangle 45,x=1.0cm,y=1.0cm]
\draw [->-] (1.04,-0.02)--(0.,0.) ;
\draw [->-] (1.9257164904935398,0.5254413794971975)--(1.04,-0.02) ;
\draw [->-,color=red] (0.5895515277177357,3.606679441322205)--(1.5817895691399835,3.294487722734828) ;
\draw [->-] (-0.4139558159781702,3.332867441679409)-- (0.5895515277177357,3.606679441322205);
\draw [->-] (-1.1101259707683522,2.559985099029886)--(-0.4139558159781702,3.332867441679409) ;
\draw [->-] (-1.2779298314101877,1.5334170962391824)--(-1.1101259707683522,2.559985099029886) ;
\draw [->-] (-0.8640908577953159,0.579091520810445)--(-1.2779298314101877,1.5334170962391824) ;
\draw [->-,color=red] (0.,0.)--(-0.8640908577953159,0.579091520810445) ;

\draw [->-] (0.5895515277177357,3.606679441322205)-- (-1.2779298314101877,1.5334170962391824);
\draw [->-,color=red] (1.5817895691399835,3.294487722734828)-- (-1.1101259707683522,2.559985099029886);
\draw [->-] (1.9257164904935398,0.5254413794971975)-- (1.5817895691399835,3.294487722734828);
\draw [->-,color=red] (-0.8640908577953159,0.579091520810445)-- (1.9257164904935398,0.5254413794971975);
\draw [->-] (-1.2779298314101877,1.5334170962391824)-- (1.04,-0.02);
\draw [->-,color=red] (-1.1101259707683522,2.559985099029886)-- (0.,0.);
\draw [->-] (-0.4139558159781702,3.332867441679409)-- (-0.8640908577953159,0.579091520810445);
\draw (0.35,4.22) node[anchor=north west] {$v_0$};
\draw (1.75,3.7) node[anchor=north west] {$v_1$};
\draw (1.85,0.5) node[anchor=north west] {$v_4$};
\draw (.85,-0.15) node[anchor=north west] {$v_5$};
\draw (-.35,-0.15) node[anchor=north west] {$v_6$};
\draw (-1.4,0.5) node[anchor=north west] {$v_7$};
\draw (-2,1.65) node[anchor=north west] {$v_8$};
\draw (-1.85,2.75) node[anchor=north west] {$v_9$};
\draw (-1.4,3.7) node[anchor=north west] {$v_{10}$};
\begin{scriptsize}
\draw [fill=black] (1.04,-0.02) circle (2.5pt);
\draw [fill=black] (1.9257164904935398,0.5254413794971975) circle (2.0pt);
\draw [fill=black] (1.5817895691399835,3.294487722734828) circle (2.0pt);
\draw [fill=black] (0.5895515277177357,3.606679441322205) circle (2.0pt);
\draw [fill=black] (-0.4139558159781702,3.332867441679409) circle (2.0pt);
\draw [fill=black] (-1.1101259707683522,2.559985099029886) circle (2.0pt);
\draw [fill=black] (-1.2779298314101877,1.5334170962391824) circle (2.0pt);
\draw [fill=black] (-0.8640908577953159,0.579091520810445) circle (2.0pt);
\draw [fill=black] (0.,0.) circle (2.0pt);
\end{scriptsize}
\end{tikzpicture}
\caption{Subdigraph 2 of $\Gamma(\mathbb{Z}_{11},\{1,8\})$}
\label{Figure: Circ11 Case2}
\end{figure}

\textbf{Case 2:} Lastly, suppose that $v_{i-1+4(3^{-1})}$ is not operational. Then, $v_{i+1}$ must not be operational since both of its in-neighbours, $v_{i-1+4(3^{-1})}$ and $v_i$, are not operational. Similarly, the vertex $v_{i-2+4(3^{-1})}$ must not be operational since both of its out-neighbours, $v_{i-1+4(3^{-1})}$ and $v_i$, are not operational. Also, since each operational vertex has an operational out-neighbour, the vertex $v_{i+2(3^{-1})}$ is operational. 

This accounts for four of the five non-operational vertices. Note that if $v_{i+4(3^{-1})}$ is not operational, then $v_{i+2}$ is operational and has no operational in-neighbours. Therefore, $v_{i+4(3^{-1})}$ must be operational. Thus, if $v_{i+2}$ is operational then the following directed path exists. $$\{v_0,v_1,\dots,v_{i-1},v_{i-1+2(3^{-1})},v_{i+2(3^{-1})},v_{i+4(3^{-1})},v_{i+2}\}$$ For illustration, the case where $v_2, v_3, v_5$, and $v_6$ are not operational for $\Gamma(\mathbb{Z}_{11},\{1,8\})$  is given in Figure \ref{Figure: Circ11 Case3}.

\begin{figure}[h!]
\centering
\begin{tikzpicture}[scale=.95, line cap=round,line join=round,>=triangle 45,x=1.0cm,y=1.0cm]
\draw [->-,color=red] (0.5895515277177357,3.606679441322205)--(1.5817895691399835,3.294487722734828) ;
\draw [->-] (-0.4139558159781702,3.332867441679409)-- (0.5895515277177357,3.606679441322205);
\draw [->-,color=red] (-1.1101259707683522,2.559985099029886)--(-0.4139558159781702,3.332867441679409) ;
\draw [->-] (-1.2779298314101877,1.5334170962391824)--(-1.1101259707683522,2.559985099029886) ;
\draw [->-] (-0.8640908577953159,0.579091520810445)--(-1.2779298314101877,1.5334170962391824) ;
\draw [->-] (0.5895515277177357,3.606679441322205)-- (-1.2779298314101877,1.5334170962391824);
\draw [->-,color=red] (1.5817895691399835,3.294487722734828)-- (-1.1101259707683522,2.559985099029886);
\draw [->-] (1.9257164904935398,0.5254413794971975)-- (1.5817895691399835,3.294487722734828);
\draw [->-,color=red] (-0.8640908577953159,0.579091520810445)-- (1.9257164904935398,0.5254413794971975);
\draw [->-,color=red] (-0.4139558159781702,3.332867441679409)-- (-0.8640908577953159,0.579091520810445);
\draw (0.35,4.22) node[anchor=north west] {$v_0$};
\draw (1.75,3.7) node[anchor=north west] {$v_1$};
\draw (1.85,0.5) node[anchor=north west] {$v_4$};
\draw (-1.4,0.5) node[anchor=north west] {$v_7$};
\draw (-2,1.65) node[anchor=north west] {$v_8$};
\draw (-1.85,2.75) node[anchor=north west] {$v_9$};
\draw (-1.4,3.7) node[anchor=north west] {$v_{10}$};
\begin{scriptsize}
\draw [fill=black] (1.9257164904935398,0.5254413794971975) circle (2.0pt);
\draw [fill=black] (1.5817895691399835,3.294487722734828) circle (2.0pt);
\draw [fill=black] (0.5895515277177357,3.606679441322205) circle (2.0pt);
\draw [fill=black] (-0.4139558159781702,3.332867441679409) circle (2.0pt);
\draw [fill=black] (-1.1101259707683522,2.559985099029886) circle (2.0pt);
\draw [fill=black] (-1.2779298314101877,1.5334170962391824) circle (2.0pt);
\draw [fill=black] (-0.8640908577953159,0.579091520810445) circle (2.0pt);

\end{scriptsize}
\end{tikzpicture}
\caption{Subdigraph 3 of $\Gamma(\mathbb{Z}_{11},\{1,8\})$}
\label{Figure: Circ11 Case3}
\end{figure}

If $v_{i+2}$ is not operational, then all of the non-operational vertices are determined. Note, when $n=11$, $-2+4(3^{-1})\equiv 3$ (mod $11$) and $-1+4(3^{-1})\equiv i+4$(mod $11$). Thus, the five non operational vertices $\{v_i,v_{i+1},v_{i+2},v_{i-1+4(3^{-1})},v_{i-2+4(3^{-1})}\}$ are in fact $\{v_i,v_{i+1},v_{i+2},v_{i+3},v_{i+4}\}$. However, as mentioned in the first case, this is a strongly connected subdigraph, as illustrated in Figure \ref{Figure: Circ11 Case4}.

\begin{figure}[h!]
\centering
\begin{tikzpicture}[scale=.95, line cap=round,line join=round,>=triangle 45,x=1.0cm,y=1.0cm]
\draw [->-] (0.5895515277177357,3.606679441322205)--(1.5817895691399835,3.294487722734828) ;
\draw [->-] (-0.4139558159781702,3.332867441679409)-- (0.5895515277177357,3.606679441322205);
\draw [->-] (-1.1101259707683522,2.559985099029886)--(-0.4139558159781702,3.332867441679409) ;
\draw [->-] (-1.2779298314101877,1.5334170962391824)--(-1.1101259707683522,2.559985099029886) ;
\draw [->-] (-0.8640908577953159,0.579091520810445)--(-1.2779298314101877,1.5334170962391824) ;
\draw [->-] (0.5895515277177357,3.606679441322205)-- (-1.2779298314101877,1.5334170962391824);
\draw [->-] (1.5817895691399835,3.294487722734828)-- (-1.1101259707683522,2.559985099029886);
\draw [->-] (-0.4139558159781702,3.332867441679409)-- (-0.8640908577953159,0.579091520810445);
\draw (0.35,4.22) node[anchor=north west] {$v_0$};
\draw (1.75,3.7) node[anchor=north west] {$v_1$};
\draw (-1.4,0.5) node[anchor=north west] {$v_7$};
\draw (-2,1.65) node[anchor=north west] {$v_8$};
\draw (-1.85,2.75) node[anchor=north west] {$v_9$};
\draw (-1.4,3.7) node[anchor=north west] {$v_{10}$};
\begin{scriptsize}
\draw [fill=black] (1.5817895691399835,3.294487722734828) circle (2.0pt);
\draw [fill=black] (0.5895515277177357,3.606679441322205) circle (2.0pt);
\draw [fill=black] (-0.4139558159781702,3.332867441679409) circle (2.0pt);
\draw [fill=black] (-1.1101259707683522,2.559985099029886) circle (2.0pt);
\draw [fill=black] (-1.2779298314101877,1.5334170962391824) circle (2.0pt);
\draw [fill=black] (-0.8640908577953159,0.579091520810445) circle (2.0pt);
\end{scriptsize}
\end{tikzpicture}
\caption{Subdigraph 4 of $\Gamma(\mathbb{Z}_{11},\{1,8\})$}
\label{Figure: Circ11 Case4}
\end{figure}
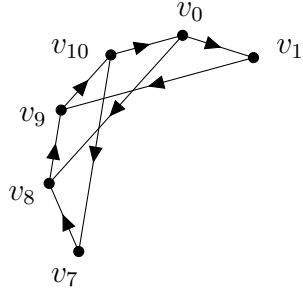
When $n>11$, the index $-2+4(3^{-1})\equiv 3$ (mod $11$) is larger than 3. Therefore, the vertices $v_{i+3}$ and $v_{i+1+4(3^{-1})}$ are both operational. Thus, the following directed path exists. $$\{v_0,v_1,\dots,v_{i-1},v_{i-1+2(3^{-1})},v_{i+2(3^{-1})},v_{i+4(3^{-1})},v_{i+1+4(3^{-1})},v_{i+3}\}$$ Therefore, if five or fewer vertices are not operational, and all of the operational vertices have at least one operational out-neighbour and one operational in-neighbour, then the only induced subdigraphs that are not strongly connected are the trivial subdigraphs.
\end{pf}

We now show that while the $F_1,F_2,F_3$ and $F_4$ for $\Gamma\left(\mathbb{Z}_{n},\{1,2(3^{-1})\}\right)$ is maximal, it may equal the $F_1$, $F_2$, $F_3$, or $F_4$ term for some other directed circulant graphs. However, we will show that the $F_5$ term for $\Gamma\left(\mathbb{Z}_{n},\{1,2(3^{-1})\}\right)$ is larger than all other circulant digraphs. Hence by Lemma \ref{Lemma: F compare}, $\Gamma\left(\mathbb{Z}_{n},\{1,2(3^{-1})\}\right)$ will be an optimally-greatest circulant graph near one.

\begin{lem}\label{Lemma: F2 Odd Part 1}
If $G=\Gamma\left(\mathbb{Z}_{n},\{1,2(3^{-1})\}\right)$ and $H= \Gamma\left(\mathbb{Z}_{n},\{a,b\}\right)$, then $F_2(G)\geq F_2(H)$.
\end{lem}
\begin{pf}
Let $G=\Gamma\left(\mathbb{Z}_{n},\{1,2(3^{-1})\}\right)$ and $H= \Gamma\left(\mathbb{Z}_{n},\{a,b\}\right)$. Suppose that only two vertices of $G$ are not operational. By Lemma \ref{Lemma: Property Odd Not Divisible By Three}, if all of the  operational vertices have at least one operational in-neighbour and one operational out-neighbour, then the resulting induced subdigraph of $G$ is strongly connected. Thus by Lemma \ref{Lemma: Circulant Out In}, it remains to count the ways to remove two out-neighbours of a vertex.

Since $n$ is odd, the additional property that two out-neighbours of a vertex are also out-neighbours of another vertex does not occur. Therefore, there are exactly $n$ ways to remove two out-neighbours of a vertex. Thus, the $F_2$ term for $G$ is maximal, in particular, it is given by $F_2(G)=\frac{n(n-3)}{2}$. However, it may be the case that $H$ obtains this maximum as well, therefore we get that $F_2(G)\geq F_2(H)$.
\end{pf}

\begin{lem}\label{Lemma: F3 Odd Part 1}
If $G=\Gamma\left(\mathbb{Z}_{n},\{1,2(3^{-1})\}\right)$ and $H= \Gamma\left(\mathbb{Z}_{n},\{a,b\}\right)$, then $F_3(G)\geq F_3(H)$.
\end{lem}
\begin{pf}
Let $G=\Gamma\left(\mathbb{Z}_{n},\{1,2(3^{-1})\}\right)$ and $H= \Gamma\left(\mathbb{Z}_{n},\{a,b\}\right)$. Suppose that only three vertices of $G$ are not operational. By Lemma \ref{Lemma: Property Odd Not Divisible By Three}, if all of the  operational vertices have at least one operational in-neighbour and one operational out-neighbour, then the resulting induced subdigraph of $G$ is strongly connected. Thus by Lemma \ref{Lemma: Circulant Out In}, it remains to count the ways to remove both out-neighbours of a vertex and one other vertex. 

Note that if we remove a vertex of $G$ and its two out neighbours, the induced subdigraph of $G$ will not be strongly connected by Lemma \ref{Lemma: Circulant Out In}. Thus, for each of the $n$ vertices of $G$, there is one way for the two out neighbours to be non-operational and then $n-3$ ways to pick the other non-operational vertex. Therefore, $G$ has at most $n(n-2)$ induced subdigraphs that are not strongly connected.

Note however, that some of these are being double counted. In particular, the vertices $v_{i+a},v_{i+b},v_{i+2a-b}$ can be chosen as the non-operational vertices for $v_i$ as well as $v_{i+a-b}$. This is because $v_{i+2a-b}$ and $v_{i+a}$ are both out-neighbours of $v_{i+a-b}$. Similarly, the vertices $v_{i+a},v_{i+b}, v_{i+2b-a}$ can be chosen as the non-operational vertices for $v_i$ as well as $v_{i+b-a}$. This is because $v_{i+2b-a}$ and $v_{i+b}$ are both out-neighbours of $v_{i+a-b}$.

Therefore, $n$ induced subdigraphs are being double counted. Hence, there are $n(n-4)$ induced subdigraphs of $G$ that are not strongly connected. Note that the counting arguments applied here apply to the other circulant graphs. Thus, there are at most $n \choose 3$ $-n(n-4)$ strongly connected subdigraphs. Since $G$ meets this maximum, it follows that $F_3(G)= \frac{n(n-1)(n-2)}{6}-n(n-4)$. However, it may be the case that $H$ obtains this maximum as well, therefore we get that $F_3(G)\geq F_3(H)$.
\end{pf}

\begin{lem}\label{Lemma: F4 Odd Part 1}
If $G=\Gamma\left(\mathbb{Z}_{n},\{1,2(3^{-1})\}\right)$ and $H= \Gamma\left(\mathbb{Z}_{n},\{a,b\}\right)$, then $F_4(G)\geq F_4(H)$.
\end{lem}
\begin{pf}
Let $G=\Gamma\left(\mathbb{Z}_{n},\{1,2(3^{-1})\}\right)$ and $H= \Gamma\left(\mathbb{Z}_{n},\{a,b\}\right)$. Suppose that four vertices of $G$ are not operational. By Lemma \ref{Lemma: Property Odd Not Divisible By Three}, if all of the  operational vertices have at least one operational in-neighbour and one operational out-neighbour, then the resulting induced subdigraph of $G$ is strongly connected. Thus, it remains to count the number of ways that both out-neighbours of a vertex and two other vertices can be non-operational. 

By Lemma \ref{Lemma: Property Odd Not Divisible By Three}, the two out-neighbours of a vertex are also the in-neighbours of another vertex. So, consider a vertex $v_i$ and its two out-neighbours $v_{i+a}$ and $v_{i+b}$. These two vertices are the in-neigbours of the vertex $v_{i+a+b}$. Therefore, if the two other non-operational vertices are not $v_{i}$ and $v_{i+a+b}$, the resulting induced subdigraph will not be strongly connected. Therefore, it is sufficient to consider the case when $v_i,v_{i+a},v_{i+b}$, and $v_{i+a+b}$ are all not operational. 

The resulting induced subdigraph of $G$ will not be strongly connected if some operational vertex has both of its out-neighbours or in-neighbours non-operational. For instance, consider $v_i$ and $v_{i+a}$. The two operational in-neighbours of $v_i$ are $v_{i-a}$ and $v_{i-b}$. The operational in-neighbour of $v_{i+a}$ is $v_{i+a-b}$. Thus, if $a-b=-a$ or $a-b=-b$, then they would share an in-neighbours. However, for $G$, $1-2(3)^{-1}\neq -1$ and $1-2(3)^{-1}\neq -2(3)^{-1}$. 

Similar calculations show that the sets $\{v_i,v_{i+b}\}$, $\{v_{i+a+b},v_{i+a}\}$, and $\{v_{i+a+b},v_{i+b}\}$ do not have any common operational out or in-neighbours. Therefore, when $v_i,v_{i+a},v_{i+b}$, and $v_{i+a+b}$ are all not operational, the resulting induced subdigraph of $G$ is strongly connected. Thus,$G$ achieves the maximum value for $F_4$ possible. However, it may be the case that $H$ obtains this maximum as well, therefore we get that $F_4(G)\geq F_4(H)$.
\end{pf}

We have showed that the $F_1$, $F_2$, $F_3$, and $F_4$ terms in the $F$-form of the strongly connected node reliability polynomial for $\Gamma\left(\mathbb{Z}_{n},\{1,2(3^{-1})\}\right)$ are the largest possible. However, there may be other directed circulant graphs that also have this property. The following theorem shows that  $\Gamma\left(\mathbb{Z}_{n},\{1,2(3^{-1})\}\right)$ has the largest $F_5$ term however. This then proves that $\Gamma\left(\mathbb{Z}_{n},\{1,2(3^{-1})\}\right)$ is an optimally-greatest circulant digraph near one.

\begin{thrm}
\label{Theorem: Odd Part 1 Circulant}
If $G=\Gamma\left(\mathbb{Z}_{n},\{1,2(3^{-1})\}\right)$ and $n$ is not divisible by three, then $G$ is an optimally-greatest directed circulant graph of the form $\Gamma\left(\mathbb{Z}_{n},\{a,b\}\right)$ for values of $p$ near one.
\end{thrm}
\begin{pf}
By Lemma \ref{Lemma: F2 Odd Part 1}, Lemma \ref{Lemma: F3 Odd Part 1}, and Lemma \ref{Lemma: F4 Odd Part 1}, the $F_i$ terms for $G$ are maximal for $i=2,3,4$. We will now show that the $F_5$ term for $G$ is larger than that for any other circulant digraph. Suppose that five vertices are not operational. By Lemma \ref{Lemma: Property Odd Not Divisible By Three}, if all of the operational vertices have at least one operational in-neighbour and one operational out-neighbour, then the resulting induced subdigraph of $G$ is strongly connected. Thus by Lemma \ref{Lemma: Circulant Out In}, it remains to count the ways to remove both out-neighbours and three other vertices.

So, consider a vertex $v_i$ and its two out-neighbours $v_{i+a}$ and $v_{i+b}$. These two vertices are the in-neigbours of the vertex $v_{i+a+b}$. Therefore, if the two other non-operational vertices are not $v_{i}$ and $v_{i+a+b}$, the resulting induced subdigraph will not be strongly connected. Therefore, it is sufficient to consider the case when $v_i,v_{i+a},v_{i+b}$, and $v_{i+a+b}$ are all not operational. 

Recall from the proof of Lemma \ref{Lemma: F4 Odd Part 1} that if just $v_i,v_{i+a},v_{i+b}$, and $v_{i+a+b}$ are not operational, then each operational vertex has at least one operational in-neighbour and at least one operational out-neighbour. Therefore, by Lemma \ref{Lemma: Property Odd Not Divisible By Three}, the only way the induced subdigraph will be disconnected, is if the other non-operational vertex violates this property. Using the known four not operational vertices, it can be shown that if one of the following vertices is not operational, then the property is violated:
\begin{align*}
v_{i-a+b}\\
v_{i-b+a}\\
v_{i+2a-b}\\
v_{i+2b-a}\\
v_{i+2a}\\
v_{i+2b}
\end{align*}

For instance, if $v_{i-a+b}$ is not operational, then the vertex $v_{i-a}$ has no operational out-neighbours since $v_{i}$ and $v_{i-a+b}$ are not operational. Now, for particular $a$ and $b$, some of these vertices indices are equal. Then, since all of the other counting was equivalent, if $G$ has the most of these indices being the same, then $G$ will have the fewest induced subdigraphs that are not strongly connected. Hence, $G$ will have the largest $F_5$ term among all directed circulant graphs. By equating every pair of indices, the following information is obtained.
\begin{align*}
&\text{If $a=2(3)^{-1}b$ then $v_{i-a+b}=v_{i+2a-b}$ and $v_{i+2b-a}=v_{i+2a}$}\\
&\text{If $a=3(2)^{-1}b$ then $v_{i-b+a}=v_{i+2b-a}$ and $v_{i+2a-b}=v_{i+2b}$}\\
&\text{If $a=3b$ then $v_{i-b+a}=v_{i+2b}$}\\
&\text{If $a=3^{-1}b$ then $v_{i-a+b}=v_{i+2a}$}\\
&\text{If $a=-b$ then $v_{i-a+b}=v_{i+2b}$ and $v_{i-b+a}=v_{i+2a}$}
\end{align*}

Therefore, up to isomorphism, the two cases to consider are when $a=2(3)^{-1}b$ or when $a=-b$. When $a=-b$ however, the circulant digraph $\Gamma(\mathbb{Z}_n,\{a,-a\}$ is the one whose underlying graph is a cycle with edges replaced by bundles. This digraph has $F_2=n$, which is lower than every other connected circulant digraph. Therefore, it is sufficient to consider directed circulant graphs of the form $\Gamma(\mathbb{Z}_n,\{b,2(3)^{-1}b\}$.

If $b$ is a unit, then $b^{-1}\{b,2(3)^{-1}b\}=\{1,2(3)^{-1}\}$. Therefore, by Lemma \ref{Lemma: Circ Iso}, it follows that $\Gamma(\mathbb{Z}_n,\{b,2(3)^{-1}b\}\cong \Gamma(\mathbb{Z}_n,\{1,2(3)^{-1}\}$. If $b$ is not a unit, then $b$ and $n$ share a common factor $r>1$. Therefore, it following that $\gcd(n,b,2(3)^{-1}b)=r$. Thus, by Lemma \ref{Lemma: Circ Disconnected}, these circulant graphs are disconnected. In summary, up to isomorphism, $\Gamma(\mathbb{Z}_n,\{1,2(3)^{-1}\}$ has the largest $F_5$ term. Therefore, by Lemma \ref{Lemma: F compare}, it is an optimally-greatest circulant digraph of the form $\Gamma(\mathbb{Z}_n,\{a,b\}$ near one when $n$ is odd and not divisible by three.
\end{pf}

\subsection{The order of the graph is odd and divisible by three}

In this section, we will show that if $n$ is odd and divisible by three, then $\Gamma(\mathbb{Z}_n,\{1,3(2)^{-1}\}$ is optimally-greatest for values of $p$ sufficiently close to one. We will first show that if five or fewer vertices are not operational, having an operational vertex with no operational out-neighbours or no operational in-neighbours is the only way for  $\Gamma(\mathbb{Z}_n,\{1,3(2)^{-1}\}$ to have an induced subdigraph that is not strongly connected. 

\begin{lem}
\label{Lemma: Property Odd Divisible By Three}
If five or fewer vertices are not operational and all of the  operational vertices have at least one operational in-neighbour and one operational out-neighbour, then the resulting induced subdigraph of $\Gamma\left(\mathbb{Z}_{n},\{1,3(2^{-1})\}\right)$ is strongly connected.
\end{lem}
\begin{pf}
Since $n\equiv 0(\textrm{mod}~3)$ and $3(2)^{-1}$ is divisible by three, $\Gamma(\mathbb{Z}_n,\{1,3(2)^{-1}\})$ can be described as 3 copies of $C_m$, where $m=\frac{n}{3}$. The vertex $v_i$ is mapped to $v_{m,r}$, where $i=nm+r$. For example, this description of $\Gamma\left(\mathbb{Z}_{21},\{1,12\}\right)$ is shown in Figure \ref{Figure: Circ21}. We use this alternative description of $\Gamma\left(\mathbb{Z}_{n},\{1,3(2^{-1})\}\right)$ to create the directed paths between the operational vertices.

\begin{figure}[h!]
\centering
\begin{tikzpicture}[line cap=round,line join=round,>=triangle 45,x=1.0cm,y=1.0cm]
\draw [->-]  (0.,0.)-- (0.,2.);
\draw [->-] (0.,2.)-- (0.,4.);
\draw [->-] (0.,4)-- (1.,0.);
\draw [->-] (1,4)-- (2,0.);
\draw [->-] (2,4)-- (3,0.);
\draw [->-] (3,4)-- (4,0.);
\draw [->-] (4,4)-- (5,0.);
\draw [->-] (5,4)-- (6,0.);
\draw [->-] (6,4)-- (0,0.);
\draw [-<-] (1.,4.)-- (1.,2.);
\draw [-<-] (1.,2.)-- (1.,0.);
\draw [->-] (2.,0.)-- (2.,2.);
\draw [->-] (2.,2.)-- (2.,4.);
\draw [-<-] (3.,4.)-- (3.,2.);
\draw [-<-] (3.,2.)-- (3.,0.);
\draw [->-] (4.,0.)-- (4.,2.);
\draw [->-] (4.,2.)-- (4.,4.);
\draw [-<-] (5.,4.)-- (5.,2.);
\draw [-<-] (5.,2.)-- (5.,0.);
\draw [->-] (6.,0.)-- (6.,2.);
\draw [->-] (6.,2.)-- (6.,4.);
\draw [->-] (0.,0.)-- (1.,0.);
\draw [->-] (1.,0.)-- (2.,0.);
\draw [->-] (2.,0.)-- (3.,0.);
\draw [->-] (3.,0.)-- (4.,0.);
\draw [->-] (4.,0.)-- (5.,0.);
\draw [->-] (5.,0.)-- (6.,0.);
\draw [->-] (0.,2.)-- (1.,2.);
\draw [->-] (1.,2.)-- (2.,2.);
\draw [->-] (2.,2.)-- (3.,2.);
\draw [->-] (3.,2.)-- (4.,2.);
\draw [->-] (4.,2.)-- (5.,2.);
\draw [->-] (5.,2.)-- (6.,2.);
\draw [->-] (0.,4.)-- (1.,4.);
\draw [->-] (1.,4.)-- (2.,4.);
\draw [-<-] (3.,4.)-- (2.,4.);
\draw [->-] (3.,4.)-- (4.,4.);
\draw [->-] (4.,4.)-- (5.,4.);
\draw [->-] (5.,4.)-- (6.,4.);
\draw [shift={(3.,-2.9364516129032237)},->-]  plot[domain=1.1625921622665343:1.979000491323259,variable=\t]({1.*7.557404381012553*cos(\t r)+0.*7.557404381012553*sin(\t r)},{0.*7.557404381012553*cos(\t r)+1.*7.557404381012553*sin(\t r)});
\draw [shift={(3.,-7.135)},->-]  plot[domain=1.2534858024220938:1.8881068511676995,variable=\t]({1.*9.615*cos(\t r)+0.*9.615*sin(\t r)},{0.*9.615*cos(\t r)+1.*9.615*sin(\t r)});
\draw [shift={(3.,4.798837209302326)},-<-]  plot[domain=4.153680770538552:5.271097190230827,variable=\t]({1.*5.659402668249057*cos(\t r)+0.*5.659402668249057*sin(\t r)},{0.*5.659402668249057*cos(\t r)+1.*5.659402668249057*sin(\t r)});
\draw (-.15,-0.15) node[anchor=north west] {$v_{0,0}$};
\draw (0.85,-0.15) node[anchor=north west] {$v_{0,1}$};
\draw (1.85,-0.15) node[anchor=north west] {$v_{0,2}$};
\draw (2.85,-0.15) node[anchor=north west] {$v_{0,3}$};
\draw (3.85,-0.15) node[anchor=north west] {$v_{0,4}$};
\draw (4.85,-0.15) node[anchor=north west] {$v_{0,5}$};
\draw (5.85,-0.15) node[anchor=north west] {$v_{0,6}$};
\draw (-0.15,1.85) node[anchor=north west] {$v_{1,0}$};
\draw (.85,1.85) node[anchor=north west] {$v_{1,1}$};
\draw (1.85,1.85) node[anchor=north west] {$v_{1,2}$};
\draw (2.85,1.85) node[anchor=north west] {$v_{1,3}$};
\draw (3.85,1.85) node[anchor=north west] {$v_{1,4}$};
\draw (4.85,1.85) node[anchor=north west] {$v_{1,5}$};
\draw (5.85,1.85) node[anchor=north west] {$v_{1,6}$};
\draw (-0.15,4.5) node[anchor=north west] {$v_{2,0}$};
\draw (0.85,4.5) node[anchor=north west] {$v_{2,1}$};
\draw (1.85,4.5) node[anchor=north west] {$v_{2,2}$};
\draw (2.85,4.5) node[anchor=north west] {$v_{2,3}$};
\draw (3.85,4.5) node[anchor=north west] {$v_{2,4}$};
\draw (4.85,4.5) node[anchor=north west] {$v_{2,5}$};
\draw (5.85,4.5) node[anchor=north west] {$v_{2,6}$};
\begin{scriptsize}
\draw [fill=black] (0.,0.) circle (2.5pt);
\draw [fill=black] (1.,0.) circle (2.5pt);
\draw [fill=black] (2.,0.) circle (2.5pt);
\draw [fill=black] (3.,0.) circle (2.5pt);
\draw [fill=black] (4.,0.) circle (2.5pt);
\draw [fill=black] (5.,0.) circle (2.5pt);
\draw [fill=black] (6.,0.) circle (2.5pt);
\draw [fill=black] (0.,2.) circle (2.5pt);
\draw [fill=black] (1.,2.) circle (2.5pt);
\draw [fill=black] (2.,2.) circle (2.5pt);
\draw [fill=black] (3.,2.) circle (2.5pt);
\draw [fill=black] (4.,2.) circle (2.5pt);
\draw [fill=black] (5.,2.) circle (2.5pt);
\draw [fill=black] (6.,2.) circle (2.5pt);
\draw [fill=black] (0.,4.) circle (2.5pt);
\draw [fill=black] (1.,4.) circle (2.5pt);
\draw [fill=black] (2.,4.) circle (2.5pt);
\draw [fill=black] (3.,4.) circle (2.5pt);
\draw [fill=black] (4.,4.) circle (2.5pt);
\draw [fill=black] (5.,4.) circle (2.5pt);
\draw [fill=black] (6.,4.) circle (2.5pt);
\end{scriptsize}
\end{tikzpicture}
\caption{Digraph $\Gamma(\mathbb{Z}_{21},\{1,12\})$}
\label{Figure: Circ21}
\end{figure}

Without loss of generality, suppose that five vertices are not operational and that $v_{0,0}$ is operational. The goal is to get a directed transversal from $v_{0,0}$ to itself by first visiting all of the operational vertices in the first copy of $C_m$. Then, we can use the same arguments developed there to create directed transversals containing all of the operational vertices in the second and third copy. Lastly, since each of these copies are connected to one another, we can connect these constructed transversals to show the induced subdigraph is strongly connected.

To start, suppose that the first non-operational vertex is in the first copy and that $v_{0,i}$ is the first nonoperational vertex. Then, since all of the operational vertices have at least one operational in-neighbour and one operational out-neighbour, the vertex $v_{1,i-1}$ is operational. 

\textbf{Case 1:} Suppose that the vertex $v_{0,i+1}$ is operational. Then, since operational vertices have an operational in neighbour, the vertex $v_{2,i}$ is operational. Thus, since one of $v_{1,i-1}$'s out-neighbours are operational, one of the directed paths exists. For illustration, the first directed path is highlighted in green in Figure \ref{Figure: Circ21 Case1}. $$\{v_{0,i-1},v_{1,i-1},v_{2,i-1},v_{2,i},v_{0,i+1}\}\hspace{30pt}\{v_{0,i-1},v_{1,i-1},v_{1,i},v_{2,i},v_{0,i+1}\}$$ 

\begin{figure}[h!]
\centering
\begin{tikzpicture}[line cap=round,line join=round,>=triangle 45,x=1.0cm,y=1.0cm]
\draw [->-]  (0.,0.)-- (0.,2.);
\draw [->-] (0.,2.)-- (0.,4.);
\draw [->-] (0.,4)-- (1.,0.);
\draw [->-,color=green] (2,4)-- (3,0.);
\draw [->-] (3,4)-- (4,0.);
\draw [->-] (4,4)-- (5,0.);
\draw [->-] (5,4)-- (6,0.);
\draw [->-] (6,4)-- (0,0.);
\draw [-<-,color=green] (1.,4.)-- (1.,2.);
\draw [-<-,color=green] (1.,2.)-- (1.,0.);

\draw [->-] (2.,2.)-- (2.,4.);
\draw [-<-] (3.,4.)-- (3.,2.);
\draw [-<-] (3.,2.)-- (3.,0.);
\draw [->-] (4.,0.)-- (4.,2.);
\draw [->-] (4.,2.)-- (4.,4.);
\draw [-<-] (5.,4.)-- (5.,2.);
\draw [-<-] (5.,2.)-- (5.,0.);
\draw [->-] (6.,0.)-- (6.,2.);
\draw [->-] (6.,2.)-- (6.,4.);
\draw [->-] (0.,0.)-- (1.,0.);

\draw [->-] (3.,0.)-- (4.,0.);
\draw [->-] (4.,0.)-- (5.,0.);
\draw [->-] (5.,0.)-- (6.,0.);
\draw [->-] (0.,2.)-- (1.,2.);
\draw [->-] (1.,2.)-- (2.,2.);
\draw [->-] (2.,2.)-- (3.,2.);
\draw [->-] (3.,2.)-- (4.,2.);
\draw [->-] (4.,2.)-- (5.,2.);
\draw [->-] (5.,2.)-- (6.,2.);
\draw [->-] (0.,4.)-- (1.,4.);
\draw [->-,color=green] (1.,4.)-- (2.,4.);
\draw [-<-] (3.,4.)-- (2.,4.);
\draw [->-] (3.,4.)-- (4.,4.);
\draw [->-] (4.,4.)-- (5.,4.);
\draw [->-] (5.,4.)-- (6.,4.);
\draw [shift={(3.,-2.9364516129032237)},->-]  plot[domain=1.1625921622665343:1.979000491323259,variable=\t]({1.*7.557404381012553*cos(\t r)+0.*7.557404381012553*sin(\t r)},{0.*7.557404381012553*cos(\t r)+1.*7.557404381012553*sin(\t r)});
\draw [shift={(3.,-7.135)},->-]  plot[domain=1.2534858024220938:1.8881068511676995,variable=\t]({1.*9.615*cos(\t r)+0.*9.615*sin(\t r)},{0.*9.615*cos(\t r)+1.*9.615*sin(\t r)});
\draw [shift={(3.,4.798837209302326)},-<-]  plot[domain=4.153680770538552:5.271097190230827,variable=\t]({1.*5.659402668249057*cos(\t r)+0.*5.659402668249057*sin(\t r)},{0.*5.659402668249057*cos(\t r)+1.*5.659402668249057*sin(\t r)});
\draw (-.15,-0.15) node[anchor=north west] {$v_{0,0}$};
\draw (0.85,-0.15) node[anchor=north west] {$v_{0,1}$};
\draw (2.85,-0.15) node[anchor=north west] {$v_{0,3}$};
\draw (3.85,-0.15) node[anchor=north west] {$v_{0,4}$};
\draw (4.85,-0.15) node[anchor=north west] {$v_{0,5}$};
\draw (5.85,-0.15) node[anchor=north west] {$v_{0,6}$};
\draw (-0.15,1.85) node[anchor=north west] {$v_{1,0}$};
\draw (.85,1.85) node[anchor=north west] {$v_{1,1}$};
\draw (1.85,1.85) node[anchor=north west] {$v_{1,2}$};
\draw (2.85,1.85) node[anchor=north west] {$v_{1,3}$};
\draw (3.85,1.85) node[anchor=north west] {$v_{1,4}$};
\draw (4.85,1.85) node[anchor=north west] {$v_{1,5}$};
\draw (5.85,1.85) node[anchor=north west] {$v_{1,6}$};
\draw (-0.15,4.5) node[anchor=north west] {$v_{2,0}$};
\draw (0.85,4.5) node[anchor=north west] {$v_{2,1}$};
\draw (1.85,4.5) node[anchor=north west] {$v_{2,2}$};
\draw (2.85,4.5) node[anchor=north west] {$v_{2,3}$};
\draw (3.85,4.5) node[anchor=north west] {$v_{2,4}$};
\draw (4.85,4.5) node[anchor=north west] {$v_{2,5}$};
\draw (5.85,4.5) node[anchor=north west] {$v_{2,6}$};
\begin{scriptsize}
\draw [fill=black] (0.,0.) circle (2.5pt);
\draw [fill=black] (1.,0.) circle (2.5pt);

\draw [fill=black] (3.,0.) circle (2.5pt);
\draw [fill=black] (4.,0.) circle (2.5pt);
\draw [fill=black] (5.,0.) circle (2.5pt);
\draw [fill=black] (6.,0.) circle (2.5pt);
\draw [fill=black] (0.,2.) circle (2.5pt);
\draw [fill=black] (1.,2.) circle (2.5pt);
\draw [fill=black] (2.,2.) circle (2.5pt);
\draw [fill=black] (3.,2.) circle (2.5pt);
\draw [fill=black] (4.,2.) circle (2.5pt);
\draw [fill=black] (5.,2.) circle (2.5pt);
\draw [fill=black] (6.,2.) circle (2.5pt);
\draw [fill=black] (0.,4.) circle (2.5pt);
\draw [fill=black] (1.,4.) circle (2.5pt);
\draw [fill=black] (2.,4.) circle (2.5pt);
\draw [fill=black] (3.,4.) circle (2.5pt);
\draw [fill=black] (4.,4.) circle (2.5pt);
\draw [fill=black] (5.,4.) circle (2.5pt);
\draw [fill=black] (6.,4.) circle (2.5pt);
\end{scriptsize}
\end{tikzpicture}
\caption{Subdigraph 1 of $\Gamma(\mathbb{Z}_{21},\{1,12\})$}
\label{Figure: Circ21 Case1}
\end{figure}

\textbf{Case 2:} Suppose that $v_{0,i+1}$ is not operational and that $v_{2,i}$ is operational. Then, since $v_{2,i}$ has an operational out-neighbour, the vertex $v_{2,i+1}$ is operational. Using the same reasoning, if $v_{0,i+2}$ is not operational, then $v_{2,i+2}$ is operational. Therefore, since there are only five not operational vertices, eventually the vertex $v_{0,i+c}$ will be operational, and a directed path from $v_{0,0}$ to $v_{0,i+c}$ will be obtained. For illustration, the case where $v_{0,2}$ and $v_{0,3}$ are not operational is given in Figure \ref{Figure: Circ21 Case2}. Again, it is assumed the first of the two directed paths highlighted in green is used.

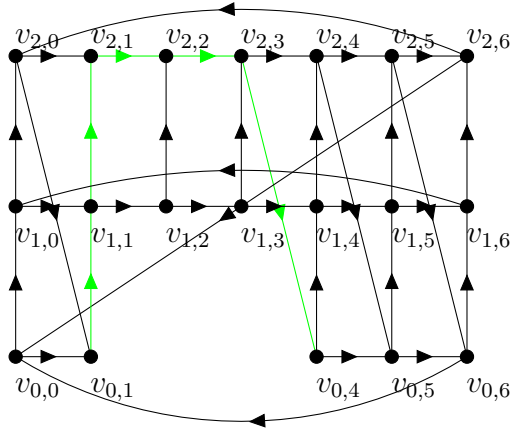
\begin{figure}[h!]
\centering
\begin{tikzpicture}[line cap=round,line join=round,>=triangle 45,x=1.0cm,y=1.0cm]
\draw [->-]  (0.,0.)-- (0.,2.);
\draw [->-] (0.,2.)-- (0.,4.);
\draw [->-] (0.,4)-- (1.,0.);

\draw [->-,color=green] (3,4)-- (4,0.);
\draw [->-] (4,4)-- (5,0.);
\draw [->-] (5,4)-- (6,0.);
\draw [->-] (6,4)-- (0,0.);
\draw [-<-,color=green] (1.,4.)-- (1.,2.);
\draw [-<-,color=green] (1.,2.)-- (1.,0.);

\draw [->-] (2.,2.)-- (2.,4.);
\draw [-<-] (3.,4.)-- (3.,2.);

\draw [->-] (4.,0.)-- (4.,2.);
\draw [->-] (4.,2.)-- (4.,4.);
\draw [-<-] (5.,4.)-- (5.,2.);
\draw [-<-] (5.,2.)-- (5.,0.);
\draw [->-] (6.,0.)-- (6.,2.);
\draw [->-] (6.,2.)-- (6.,4.);
\draw [->-] (0.,0.)-- (1.,0.);

\draw [->-] (4.,0.)-- (5.,0.);
\draw [->-] (5.,0.)-- (6.,0.);
\draw [->-] (0.,2.)-- (1.,2.);
\draw [->-] (1.,2.)-- (2.,2.);
\draw [->-] (2.,2.)-- (3.,2.);
\draw [->-] (3.,2.)-- (4.,2.);
\draw [->-] (4.,2.)-- (5.,2.);
\draw [->-] (5.,2.)-- (6.,2.);
\draw [->-] (0.,4.)-- (1.,4.);
\draw [->-,color=green] (1.,4.)-- (2.,4.);
\draw [-<-,,color=green] (3.,4.)-- (2.,4.);
\draw [->-] (3.,4.)-- (4.,4.);
\draw [->-] (4.,4.)-- (5.,4.);
\draw [->-] (5.,4.)-- (6.,4.);
\draw [shift={(3.,-2.9364516129032237)},->-]  plot[domain=1.1625921622665343:1.979000491323259,variable=\t]({1.*7.557404381012553*cos(\t r)+0.*7.557404381012553*sin(\t r)},{0.*7.557404381012553*cos(\t r)+1.*7.557404381012553*sin(\t r)});
\draw [shift={(3.,-7.135)},->-]  plot[domain=1.2534858024220938:1.8881068511676995,variable=\t]({1.*9.615*cos(\t r)+0.*9.615*sin(\t r)},{0.*9.615*cos(\t r)+1.*9.615*sin(\t r)});
\draw [shift={(3.,4.798837209302326)},-<-]  plot[domain=4.153680770538552:5.271097190230827,variable=\t]({1.*5.659402668249057*cos(\t r)+0.*5.659402668249057*sin(\t r)},{0.*5.659402668249057*cos(\t r)+1.*5.659402668249057*sin(\t r)});
\draw (-.15,-0.15) node[anchor=north west] {$v_{0,0}$};
\draw (0.85,-0.15) node[anchor=north west] {$v_{0,1}$};
\draw (3.85,-0.15) node[anchor=north west] {$v_{0,4}$};
\draw (4.85,-0.15) node[anchor=north west] {$v_{0,5}$};
\draw (5.85,-0.15) node[anchor=north west] {$v_{0,6}$};
\draw (-0.15,1.85) node[anchor=north west] {$v_{1,0}$};
\draw (.85,1.85) node[anchor=north west] {$v_{1,1}$};
\draw (1.85,1.85) node[anchor=north west] {$v_{1,2}$};
\draw (2.85,1.85) node[anchor=north west] {$v_{1,3}$};
\draw (3.85,1.85) node[anchor=north west] {$v_{1,4}$};
\draw (4.85,1.85) node[anchor=north west] {$v_{1,5}$};
\draw (5.85,1.85) node[anchor=north west] {$v_{1,6}$};
\draw (-0.15,4.5) node[anchor=north west] {$v_{2,0}$};
\draw (0.85,4.5) node[anchor=north west] {$v_{2,1}$};
\draw (1.85,4.5) node[anchor=north west] {$v_{2,2}$};
\draw (2.85,4.5) node[anchor=north west] {$v_{2,3}$};
\draw (3.85,4.5) node[anchor=north west] {$v_{2,4}$};
\draw (4.85,4.5) node[anchor=north west] {$v_{2,5}$};
\draw (5.85,4.5) node[anchor=north west] {$v_{2,6}$};
\begin{scriptsize}
\draw [fill=black] (0.,0.) circle (2.5pt);
\draw [fill=black] (1.,0.) circle (2.5pt);

\draw [fill=black] (4.,0.) circle (2.5pt);
\draw [fill=black] (5.,0.) circle (2.5pt);
\draw [fill=black] (6.,0.) circle (2.5pt);
\draw [fill=black] (0.,2.) circle (2.5pt);
\draw [fill=black] (1.,2.) circle (2.5pt);
\draw [fill=black] (2.,2.) circle (2.5pt);
\draw [fill=black] (3.,2.) circle (2.5pt);
\draw [fill=black] (4.,2.) circle (2.5pt);
\draw [fill=black] (5.,2.) circle (2.5pt);
\draw [fill=black] (6.,2.) circle (2.5pt);
\draw [fill=black] (0.,4.) circle (2.5pt);
\draw [fill=black] (1.,4.) circle (2.5pt);
\draw [fill=black] (2.,4.) circle (2.5pt);
\draw [fill=black] (3.,4.) circle (2.5pt);
\draw [fill=black] (4.,4.) circle (2.5pt);
\draw [fill=black] (5.,4.) circle (2.5pt);
\draw [fill=black] (6.,4.) circle (2.5pt);
\end{scriptsize}
\end{tikzpicture}
\caption{Subdigraph 2 of $\Gamma(\mathbb{Z}_{21},\{1,12\})$}
\label{Figure: Circ21 Case2}
\end{figure}

\textbf{Case 3:} Lastly, suppose that $v_{0,i+1}$ is not operational and $v_{2,i}$ is not operational. Then, since $v_{2,i}$ and $v_{0,i}$ are both not operational out neighbours of $v_{2,i-1}$, the vertex $v_{2,i-1}$ must be not operational. Now, notice that if $v_{2,i+1}$ is not operational, then all five of the non-operational vertices are decided and the operational vertex $v_{0,i+2}$ would have no operational in-neighbours. Therefore, the following directed path exists. $$\{v_{0,i-1},v_{1,i-1},v_{1,i},v_{1,i+2},v_{2,i+1}\}$$ For illustration, this directed path is highlighted in green in Figure \ref{Figure: Circ21 Case3}.

\begin{figure}[h!]
\centering
\begin{tikzpicture}[line cap=round,line join=round,>=triangle 45,x=1.0cm,y=1.0cm]
\draw [->-]  (0.,0.)-- (0.,2.);
\draw [->-] (0.,2.)-- (0.,4.);
\draw [->-] (0.,4)-- (1.,0.);
\draw [->-] (3,4)-- (4,0.);
\draw [->-] (4,4)-- (5,0.);
\draw [->-] (5,4)-- (6,0.);
\draw [->-] (6,4)-- (0,0.);
\draw [-<-,color=green] (1.,2.)-- (1.,0.);
\draw [-<-,color=green] (3.,4.)-- (3.,2.);
\draw [->-] (4.,0.)-- (4.,2.);
\draw [->-] (4.,2.)-- (4.,4.);
\draw [-<-] (5.,4.)-- (5.,2.);
\draw [-<-] (5.,2.)-- (5.,0.);
\draw [->-] (6.,0.)-- (6.,2.);
\draw [->-] (6.,2.)-- (6.,4.);
\draw [->-] (0.,0.)-- (1.,0.);
\draw [->-] (4.,0.)-- (5.,0.);
\draw [->-] (5.,0.)-- (6.,0.);
\draw [->-] (0.,2.)-- (1.,2.);
\draw [->-,color=green] (1.,2.)-- (2.,2.);
\draw [->-,color=green] (2.,2.)-- (3.,2.);
\draw [->-] (3.,2.)-- (4.,2.);
\draw [->-] (4.,2.)-- (5.,2.);
\draw [->-] (5.,2.)-- (6.,2.);
\draw [->-] (3.,4.)-- (4.,4.);
\draw [->-] (4.,4.)-- (5.,4.);
\draw [->-] (5.,4.)-- (6.,4.);
\draw [shift={(3.,-2.9364516129032237)},->-]  plot[domain=1.1625921622665343:1.979000491323259,variable=\t]({1.*7.557404381012553*cos(\t r)+0.*7.557404381012553*sin(\t r)},{0.*7.557404381012553*cos(\t r)+1.*7.557404381012553*sin(\t r)});
\draw [shift={(3.,-7.135)},->-]  plot[domain=1.2534858024220938:1.8881068511676995,variable=\t]({1.*9.615*cos(\t r)+0.*9.615*sin(\t r)},{0.*9.615*cos(\t r)+1.*9.615*sin(\t r)});
\draw [shift={(3.,4.798837209302326)},-<-]  plot[domain=4.153680770538552:5.271097190230827,variable=\t]({1.*5.659402668249057*cos(\t r)+0.*5.659402668249057*sin(\t r)},{0.*5.659402668249057*cos(\t r)+1.*5.659402668249057*sin(\t r)});
\draw (-.15,-0.15) node[anchor=north west] {$v_{0,0}$};
\draw (0.85,-0.15) node[anchor=north west] {$v_{0,1}$};
\draw (3.85,-0.15) node[anchor=north west] {$v_{0,4}$};
\draw (4.85,-0.15) node[anchor=north west] {$v_{0,5}$};
\draw (5.85,-0.15) node[anchor=north west] {$v_{0,6}$};
\draw (-0.15,1.85) node[anchor=north west] {$v_{1,0}$};
\draw (.85,1.85) node[anchor=north west] {$v_{1,1}$};
\draw (1.85,1.85) node[anchor=north west] {$v_{1,2}$};
\draw (2.85,1.85) node[anchor=north west] {$v_{1,3}$};
\draw (3.85,1.85) node[anchor=north west] {$v_{1,4}$};
\draw (4.85,1.85) node[anchor=north west] {$v_{1,5}$};
\draw (5.85,1.85) node[anchor=north west] {$v_{1,6}$};
\draw (-0.15,4.5) node[anchor=north west] {$v_{2,0}$};
\draw (2.85,4.5) node[anchor=north west] {$v_{2,3}$};
\draw (3.85,4.5) node[anchor=north west] {$v_{2,4}$};
\draw (4.85,4.5) node[anchor=north west] {$v_{2,5}$};
\draw (5.85,4.5) node[anchor=north west] {$v_{2,6}$};
\begin{scriptsize}
\draw [fill=black] (0.,0.) circle (2.5pt);
\draw [fill=black] (1.,0.) circle (2.5pt);
\draw [fill=black] (4.,0.) circle (2.5pt);
\draw [fill=black] (5.,0.) circle (2.5pt);
\draw [fill=black] (6.,0.) circle (2.5pt);
\draw [fill=black] (0.,2.) circle (2.5pt);
\draw [fill=black] (1.,2.) circle (2.5pt);
\draw [fill=black] (2.,2.) circle (2.5pt);
\draw [fill=black] (3.,2.) circle (2.5pt);
\draw [fill=black] (4.,2.) circle (2.5pt);
\draw [fill=black] (5.,2.) circle (2.5pt);
\draw [fill=black] (6.,2.) circle (2.5pt);
\draw [fill=black] (0.,4.) circle (2.5pt);
\draw [fill=black] (3.,4.) circle (2.5pt);
\draw [fill=black] (4.,4.) circle (2.5pt);
\draw [fill=black] (5.,4.) circle (2.5pt);
\draw [fill=black] (6.,4.) circle (2.5pt);
\end{scriptsize}
\end{tikzpicture}
\caption{Subdigraph 3 of $\Gamma(\mathbb{Z}_{21},\{1,12\})$}
\label{Figure: Circ21 Case3}
\end{figure}

\textbf{Case 4:} Lastly, if $v_{0,i+2}$ is not operational, then all five of the not operational vertices are decided. Therefore, $v_{2,i+2}$ and $v_{2,i+2}$ are operational and the following directed path exists. $$\{v_{0,i-1},v_{1,i-1},v_{1,i},v_{1,i+2},v_{2,i+1},v_{2,i+2},v_{2,i+2}\}$$ For illustration, this directed path is highlighted in green in Figure \ref{Figure: Circ21 Case4}.

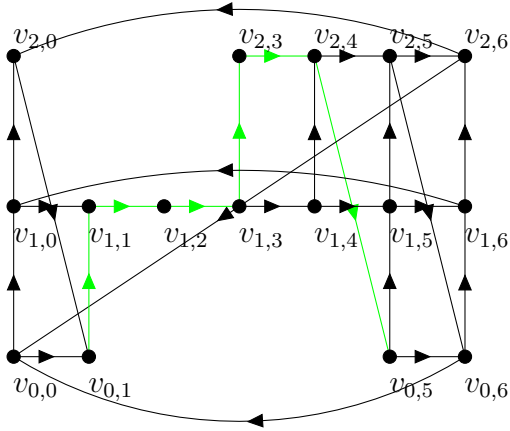
\begin{figure}[h!]
\centering
\begin{tikzpicture}[line cap=round,line join=round,>=triangle 45,x=1.0cm,y=1.0cm]
\draw [->-]  (0.,0.)-- (0.,2.);
\draw [->-] (0.,2.)-- (0.,4.);
\draw [->-] (0.,4)-- (1.,0.);
\draw [->-,color=green] (4,4)-- (5,0.);
\draw [->-] (5,4)-- (6,0.);
\draw [->-] (6,4)-- (0,0.);
\draw [-<-,color=green] (1.,2.)-- (1.,0.);
\draw [-<-,color=green] (3.,4.)-- (3.,2.);
\draw [->-] (4.,2.)-- (4.,4.);
\draw [-<-] (5.,4.)-- (5.,2.);
\draw [-<-] (5.,2.)-- (5.,0.);
\draw [->-] (6.,0.)-- (6.,2.);
\draw [->-] (6.,2.)-- (6.,4.);
\draw [->-] (0.,0.)-- (1.,0.);
\draw [->-] (5.,0.)-- (6.,0.);
\draw [->-] (0.,2.)-- (1.,2.);
\draw [->-,color=green] (1.,2.)-- (2.,2.);
\draw [->-,color=green] (2.,2.)-- (3.,2.);
\draw [->-] (3.,2.)-- (4.,2.);
\draw [->-] (4.,2.)-- (5.,2.);
\draw [->-] (5.,2.)-- (6.,2.);
\draw [->-,color=green] (3.,4.)-- (4.,4.);
\draw [->-] (4.,4.)-- (5.,4.);
\draw [->-] (5.,4.)-- (6.,4.);
\draw [shift={(3.,-2.9364516129032237)},->-]  plot[domain=1.1625921622665343:1.979000491323259,variable=\t]({1.*7.557404381012553*cos(\t r)+0.*7.557404381012553*sin(\t r)},{0.*7.557404381012553*cos(\t r)+1.*7.557404381012553*sin(\t r)});
\draw [shift={(3.,-7.135)},->-]  plot[domain=1.2534858024220938:1.8881068511676995,variable=\t]({1.*9.615*cos(\t r)+0.*9.615*sin(\t r)},{0.*9.615*cos(\t r)+1.*9.615*sin(\t r)});
\draw [shift={(3.,4.798837209302326)},-<-]  plot[domain=4.153680770538552:5.271097190230827,variable=\t]({1.*5.659402668249057*cos(\t r)+0.*5.659402668249057*sin(\t r)},{0.*5.659402668249057*cos(\t r)+1.*5.659402668249057*sin(\t r)});
\draw (-.15,-0.15) node[anchor=north west] {$v_{0,0}$};
\draw (0.85,-0.15) node[anchor=north west] {$v_{0,1}$};
\draw (4.85,-0.15) node[anchor=north west] {$v_{0,5}$};
\draw (5.85,-0.15) node[anchor=north west] {$v_{0,6}$};
\draw (-0.15,1.85) node[anchor=north west] {$v_{1,0}$};
\draw (.85,1.85) node[anchor=north west] {$v_{1,1}$};
\draw (1.85,1.85) node[anchor=north west] {$v_{1,2}$};
\draw (2.85,1.85) node[anchor=north west] {$v_{1,3}$};
\draw (3.85,1.85) node[anchor=north west] {$v_{1,4}$};
\draw (4.85,1.85) node[anchor=north west] {$v_{1,5}$};
\draw (5.85,1.85) node[anchor=north west] {$v_{1,6}$};
\draw (-0.15,4.5) node[anchor=north west] {$v_{2,0}$};
\draw (2.85,4.5) node[anchor=north west] {$v_{2,3}$};
\draw (3.85,4.5) node[anchor=north west] {$v_{2,4}$};
\draw (4.85,4.5) node[anchor=north west] {$v_{2,5}$};
\draw (5.85,4.5) node[anchor=north west] {$v_{2,6}$};
\begin{scriptsize}
\draw [fill=black] (0.,0.) circle (2.5pt);
\draw [fill=black] (1.,0.) circle (2.5pt);
\draw [fill=black] (5.,0.) circle (2.5pt);
\draw [fill=black] (6.,0.) circle (2.5pt);
\draw [fill=black] (0.,2.) circle (2.5pt);
\draw [fill=black] (1.,2.) circle (2.5pt);
\draw [fill=black] (2.,2.) circle (2.5pt);
\draw [fill=black] (3.,2.) circle (2.5pt);
\draw [fill=black] (4.,2.) circle (2.5pt);
\draw [fill=black] (5.,2.) circle (2.5pt);
\draw [fill=black] (6.,2.) circle (2.5pt);
\draw [fill=black] (0.,4.) circle (2.5pt);
\draw [fill=black] (3.,4.) circle (2.5pt);
\draw [fill=black] (4.,4.) circle (2.5pt);
\draw [fill=black] (5.,4.) circle (2.5pt);
\draw [fill=black] (6.,4.) circle (2.5pt);
\end{scriptsize}
\end{tikzpicture}
\caption{Subdigraph 4 of $\Gamma(\mathbb{Z}_{21},\{1,12\})$}
\label{Figure: Circ21 Case4}
\end{figure}

In all cases, there is a directed path to the next operational vertex in the first copy. Therefore, a directed transversal containing all of the operational vertices in the first copy of $C_m$ exists. Using the same arguments, a directed transversal containing all of the operational vertices in the second and third copy of $C_m$ also exists. Thus, it remains to show that these transversals can be connected.

Notice that since at most five vertices are not operational, there is at least one operational vertex in the first copy of $C_m$ has an arc to an operational vertex in the second copy. Similarly, there is at least one operational vertex in the second copy of $C_m$ has an arc to an operational vertex in the third copy. Also, there is at least one operational vertex in the third copy of $C_m$ has an arc to an operational vertex in the first copy. 

Hence, there is an arc from the first tranvsersal to the second, from the second to the third, and from the third to the first. Therefore, the three transversals can all be connected to give a transversal containing all of the operational vertices. Thus, if five or fewer vertices are not operational, and all of the operational vertices have at least one operational in-neighbour and one operational out-neighbour, then the resulting induced subdigraph of $\Gamma(\mathbb{Z}_n,\{1,3(2^{-1}\})$ is strongly connected.
\end{pf}

The same argument used for directed circulant graphs of odd order that were not divisible by three will be used here. In fact, Lemma \ref{Lemma: F2 Odd Part 1}, Lemma \ref{Lemma: F3 Odd Part 1}, and Lemma \ref{Lemma: F4 Odd Part 1} all have analogous versions for $\Gamma\left(\mathbb{Z}_{n},\{1,3(2^{-1})\}\right)$. That is, the $F_1,F_2,F_3$, and $F_4$ term for $\Gamma\left(\mathbb{Z}_{n},\{1,3(2^{-1})\}\right)$ is maximal, but may be equal to that of another circulant digraph. Thus, it is sufficient to show that the $F_5$ term for $\Gamma\left(\mathbb{Z}_{n},\{1,3(2^{-1})\}\right)$ is larger than the other circulant digraphs.

\begin{thrm}
\label{Theorem: Odd Part 2 Circulant}
If $G=\Gamma\left(\mathbb{Z}_{n},\{1,3(2^{-1})\}\right)$ and $n$ is divisible by three, then $G$ is an optimally-greatest directed circulant graph of the form $\Gamma\left(\mathbb{Z}_{n},\{a,b\}\right)$ for values of $p$ near one.
\end{thrm}
\begin{pf}
By the analogous versions of Lemma \ref{Lemma: F2 Odd Part 1}, Lemma \ref{Lemma: F3 Odd Part 1}, and Lemma \ref{Lemma: F4 Odd Part 1}, the $F_i$ terms for $G$ are maximal for $i=1,2,3,4$. We will now show that the $F_5$ term for $G$ is larger than that for any other directed circulant graph. Suppose that five vertices are not operational. By Lemma \ref{Lemma: Property Odd Divisible By Three}, if all of the operational vertices have at least one operational in-neighbour and one operational out-neighbour, then the resulting induced subdigraph of $G$ is strongly connected.

Thus by Lemma \ref{Lemma: Circulant Out In}, it remains to count the ways to remove both out-neighbours and three other vertices. Then, as argued in the proof of Theorem \ref{Theorem: Odd Part 1 Circulant}, the circulant digraph that has the most of the following indices from the following vertices the same, will have the largest $F_5$ term.
\begin{align*}
v_{i-a+b}\\
v_{i-b+a}\\
v_{i+2a-b}\\
v_{i+2b-a}\\
v_{i+2a}\\
v_{i+2b}
\end{align*}

Same as before, this problem is solved by equating every pair of indices. The same information is obtained, but for completeness, the information is provided again below.
\begin{align*}
&\text{If $a=2(3)^{-1}b$ then $v_{i-a+b}=v_{i+2a-b}$ and $v_{i+2b-a}=v_{i+2a}$}\\
&\text{If $a=3(2)^{-1}b$ then $v_{i-b+a}=v_{i+2b-a}$ and $v_{i+2a-b}=v_{i+2b}$}\\
&\text{If $a=3b$ then $v_{i-b+a}=v_{i+2b}$}\\
&\text{If $a=3^{-1}b$ then $v_{i-a+b}=v_{i+2a}$}\\
&\text{If $a=-b$ then $v_{i-a+b}=v_{i+2b}$ and $v_{i-b+a}=v_{i+2a}$}
\end{align*}

In this case, the case where $a=2(3)^{-1}$ does not exist because 3 is not a unit. The two cases to consider this time are when $a=3(2)^{-1}b$ or when $a=-b$. Again, when $a=-b$ the directed circulant graph $\Gamma(\mathbb{Z}_n,\{a,-a\}$ is the cycle with edges replaced by bundles. This graph has $F_2=n$, which is lower than every other connected circulant digraph. Therefore, it is sufficient to consider circulant digraphs of the form $\Gamma(\mathbb{Z}_n,\{b,3(2)^{-1}b\}$.

If $b$ is a unit, then $b^{-1}\{b,3(2)^{-1}b\}=\{1,3(2)^{-1}\}$. Therefore, by Lemma \ref{Lemma: Circ Iso}, it follows that $\Gamma(\mathbb{Z}_n,\{b,2(3)^{-1}b\}\cong \Gamma(\mathbb{Z}_n,\{1,3(2)^{-1}\}$. If $b$ is not a unit, then $b$ and $n$ share a common factor $r>1$. Therefore, it following that $\gcd(n,b,3(2)^{-1}b)=r$. Thus, by Lemma \ref{Lemma: Circ Disconnected}, these directed circulant graphs are disconnected. In summary, up to isomorphism, $\Gamma(\mathbb{Z}_n,\{1,3(2)^{-1}\}$ has the largest $F_5$ term. Therefore, by Lemma \ref{Lemma: F compare}, it is an optimally-greatest circulant digraph of the form $\Gamma(\mathbb{Z}_n,\{a,b\}$ near one when $n$ is odd and not divisible by three.
\end{pf}

\section{Conclusions and Open Problems}

In this paper, we proved that for values of $p$ sufficiently close to zero, the digraph $\Gamma(\mathbb{Z}_n,\{1,-1\}$ is an optimally-greatest directed circulant graph of the form $\Gamma(\mathbb{Z}_n,\{a,b\}$. For values of $p$ sufficiently close to one, an optimally-greatest digraph depended on the order of the graph. In particular, if $n$ is even, then it was $\Gamma(\mathbb{Z}_{n},\{1,\frac{n}{2}+1\})$. If $n$ is odd and not divisible by three then it was $\Gamma(\mathbb{Z}_{n},\{1,2(3^{-1})\})$. If $n$ is odd and divisible by three then it was $\Gamma(\mathbb{Z}_{n},\{1,3(2^{-1})\})$. Combining these results provides our main result.

\begin{thrm}
For strongly connected node reliability, there is no optimally-greatest digraph of the form $\Gamma(\mathbb{Z}_n,\{a,b\})$.
\end{thrm}

It remains an open problem to extend this result to all directed graphs with $n$ vertices and $2n$ arc, and show that the circulant digraphs described are in fact the optimally-greatest for $p$ near one. Based on the calculations done here, there are significant structural restrictions if a digraph other than a circulant was optimal near 1. For instance, the graph would need to be regular. Also, the graph would need to have the property that two out-neighbours of a vertex are also in-neighbours of another vertex. We conjecture that the circulants described in this paper are indeed optimally-greatest for $p$ sufficiently close to one for digraphs of order $n$ and size $2n$.

\bibliographystyle{amsplain}

\bibliography{SCNodeRel_References}

\end{document}